\documentclass{amsart}

\usepackage{amssymb}
\usepackage{bbm}
\usepackage{verbatim}
\usepackage{caption,subcaption}

\usepackage{enumitem}

\usepackage{tikz}
\usetikzlibrary{calc}
\usetikzlibrary{arrows,automata,shapes}
\usetikzlibrary{snakes}
\usetikzlibrary{patterns}
\tikzstyle{every state}=[minimum size=12pt,inner sep=0pt]

\DeclareMathOperator{\NSyn}{N\Syn}
\DeclareMathOperator{\Syn}{Syn}

\DeclareSymbolFont{rsfscript}{OMS}{rsfs}{m}{n}
\DeclareSymbolFontAlphabet{\mathrsfs}{rsfscript}
\DeclareMathOperator{\Fix}{Fix}
\DeclareMathOperator{\Aut}{Aut}	%automorphism group

\DeclareMathOperator{\Out}{Out}	%outer automorphism group
\DeclareMathOperator{\id}{\mathbbm{1}}

\DeclareMathOperator{\IMG}{\hbox{IMG}}

\DeclareMathOperator{\Sch}{Sch}
\DeclareMathOperator{\St}{Stab}
\DeclareMathOperator{\SSt}{St}
\DeclareSymbolFont{rsfscript}{OMS}{rsfs}{m}{n}
\DeclareSymbolFont{bbold}{U}{bbold}{m}{n}
\DeclareSymbolFontAlphabet{\mathbbold}{bbold}

\newtheorem{theorem}{Theorem}[section]
\newtheorem{proposition}[theorem]{Proposition}

\newtheorem{lemma}[theorem]{Lemma}
\newtheorem{corollary}[theorem]{Corollary}

\newtheorem{prob}{Problem}
\newtheorem{rem}[theorem]{Remark}

\newcommand{\la}{\langle}
\newcommand{\ra}{\rangle}

\newcommand{\wt}{\widetilde}

\newcommand{\inv}{^{-1}}

\newcommand{\oo}{\overline}
\newcommand{\QQ}{{{Q}}}		%stateset
\newcommand{\XX}{{{\Sigma}}}	%alphabet
\newcommand{\xcdot}{{{\cdot}}}		%transition
\newcommand{\xcirc}{{{\circ}}}	%production
\newcommand{\nucleus}{{{\aleph}}}
\newcommand{\lacroix}{\tikz[baseline=-.5ex]{\draw[->,>=latex] (0,0) -- (4ex,0); \draw[->,>=latex] (1.8ex,2ex) -- (1.8ex,-2ex);}}

\def\vlongrightarrow{\relbar\joinrel\longrightarrow}
\def\vvlongrightarrow{\relbar\joinrel\vlongrightarrow}
\def\vvvlongrightarrow{\relbar\joinrel\vvlongrightarrow}

\def\longmapright#1{\smash{\mathop{\vlongrightarrow}\limits^{#1}}}

\def\vvlongmapright#1{\smash{\mathop{\vvvlongrightarrow}\limits^{#1}}}

\newcommand{\maprightempty}{\longrightarrow}
\newcommand{\mapright}[1]{\smash{\stackrel{\text{\tiny{$#1$}}}{\vlongrightarrow}}}

\newcommand{\card}[1]{|#1|}
\newcommand{\dual}{\hbox{$\mathfrak{d}$}}

\newcommand{\lamp}{{\mathcal{L}}}		%lamplighter
\newcommand{\measure}{{\mu}}		%measure
\newcommand{\cof}{\mathrel{\propto}}	%cofinality
\newcommand{\aut}[1]{{\mathrsfs{#1}}}			%automaton
\newcommand{\gaut}[1]{\la\aut{#1}\ra}			%automaton GROUP
\newcommand{\gdaut}[1]{\gaut{\dual\aut{#1}}}		%automaton DUAL GROUP
\newcommand{\sgaut}[1]{\gaut{#1}^+}			%automaton SEMIGROUP
\newcommand{\auta}{\aut{M}}
\newcommand{\gauta}{\gaut{M}}				%automaton GROUP <M>
\newcommand{\gdauta}{\gdaut{M}}				%automaton DUAL GROUP <dM>
\newcommand{\sgauta}{\sgaut{M}}				%automaton SEMIGROUP <M>+
\newcommand{\Prune}[1]{\mathcal{B}(#1)}

\newcommand{\krit}{\hbox{\Large{$\kappa$}}}%set of critical points
\newcommand{\SinkA}{\hbox{${\mathcal S}_{a}$}}
\newcommand{\semig}[1]{\hbox{$\langle#1\rangle^+$}}	%{S(#1)}
\newcommand{\group}[1]{\hbox{$\langle#1\rangle$}}		%{G(#1)}

\newcommand{\ball}{B}		
\newcommand{\balle}{\ball^{\tiny\rightarrow\!\!}}

\newcommand{\pavet}[{6}]{
\filldraw[fill=w,line width=0pt] (#1+0,#2+0) -- (#1+1,#2+1) -- (#1+2,#2+0) -- cycle;
\filldraw[fill=w,line width=0pt] (#1+0,#2+0) -- (#1+1,#2+1) -- (#1,#2+2) -- cycle;
\filldraw[fill=w,line width=0pt] (#1+2,#2+0) -- (#1+1,#2+1) -- (#1+2,#2+2) -- cycle;
\filldraw[fill=w,line width=0pt] (#1+0,#2+2) -- (#1+1,#2+1) -- (#1+2,#2+2) -- cycle;
\draw[line width=1pt] (#1+0,#2+0) -- (#1+2,#2+0) -- (#1+2,#2+2) -- (#1+0,#2+2) -- cycle;
\node[]  (n) at (#1+1,#2+1.7) {#3};
\node[]  (s) at (#1+1,#2+.3)  {#4};
\node[]  (e) at (#1+1.5,#2+1) {#5};
\node[]  (w) at (#1+0.5,#2+1) {#6};
}

\newcommand{\cross}[6]{
\node[] (in) at (#1,#2)   {{$\begin{array}[b]{ccc}
		& #3            	&    \\
#4		& \lacroix    	& #6 \\
		& #5 			&            		
\end{array}$} };
}
\definecolor{r}{rgb}{1,0,0}
\definecolor{g}{rgb}{0,1,0}
\definecolor{b}{rgb}{0,0,1}
\definecolor{c}{rgb}{0,1,1}
\definecolor{y}{rgb}{1,1,0}
\definecolor{n}{rgb}{0,0,0}
\definecolor{w}{rgb}{1,1,1}

\title[Boundary action of autom. groups w/o singular points \& Wang tilings~~~~]
{Boundary action of automaton groups without~singular~points and Wang tilings}

\date{\today}
\begin{document}

\author[D. D'Angeli]	{Daniele D'Angeli}
\address{Institut f\"{u}r Diskrete Mathematik\\
		Technische Universit\"{a}t Graz\\
		Steyrergasse 30, 8010 Graz, Austria.}
\email{dangeli@math.tugraz.at}

\author[Th. Godin]	{Thibault Godin}%
\address{IRIF, UMR 8243 Universit\'e Paris Diderot \& CNRS\\
		B\^atiment Sophie Germain\\
		75205 Paris Cedex 13, France.}
\email{godin@irif.fr}

\author[I. Klimann]	{Ines Klimann}
\address{IRIF, UMR 8243 Universit\'e Paris Diderot \& CNRS\\
		B\^atiment Sophie Germain\\
		75205 Paris Cedex 13, France.}
\email{klimann@irif.fr}

\author[M. Picantin]	{Matthieu Picantin}
\address{IRIF, UMR 8243 Universit\'e Paris Diderot \& CNRS\\
		B\^atiment Sophie Germain\\
		75205 Paris Cedex 13, France.}
\email{picantin@irif.fr}

\author[E. Rodaro]	{Emanuele Rodaro}
\address{Department of mathematics, Politecnico di Milano\\
		Piazza Leonardo da Vinci, 32, 20133 Milano, Italy}
\email{emanuele.rodaro@polimi.it}

\maketitle

\begin{abstract}
\vspace*{-18pt}%%%%%%%%%%%%%%%%%%%%%%%%%%%%%%%%%%%
We study automaton groups without singular points, that is, points in the boundary for which the map that associates to each point its stabilizer, is not continuous. This is motivated by the problem of finding examples of infinite bireversible automaton groups with all trivial stabilizers in the boundary, raised by Grigorchuk and Savchuk.
We show that, in general, the set of singular points has measure zero. Then we focus our attention on several classes of automata. 
We characterize those contracting automata generating groups without singular points, 
and apply this characterization to the Basilica group.
We prove that potential examples of reversible automata generating infinite groups without singular points are necessarily bireversible.
Then we provide some necessary conditions for such examples to exist, and 
study some dynamical properties of their Schreier graphs in the boundary.
Finally we relate some of those automata with aperiodic tilings of the discrete plane via Wang tilings. This 
has a series of consequences from the algorithmic and dynamical points of view, and 
is related to a problem of Gromov regarding the searching for examples of~CAT(0) complexes
whose fundamental groups are not hyperbolic and contain no subgroup isomorphic to~$\mathbb{Z}^{2}$.
\end{abstract}

{\footnotesize
\tableofcontents
}

%-------------------------------------------------------------------------------------------------------------------------------------------------
%-------------------------------------------------------------------------------------------------------------------------------------------------
%-------------------------------------------------------------------------------------------------------------------------------------------------
\section{Introduction}
The motivation comes from the study of the dynamical system $(G, \partial T, \measure)$ given by the measure $\mu$ preserving action of a group $G$ on the boundary $\partial T$ of a rooted tree $T$. By considering one orbit of this action (\emph{i.e.} a Schreier graph), one may ask if it is possible to recover the information about the original dynamics in terms of the information contained in a typical orbit. This problem may be rephrased as follows: what  conditions have to be imposed on the dynamical system $(G,\partial T, \measure)$ in order to guarantee that, for a typical point $\xi\in \partial T$, this dynamical system is isomorphic to the system $(G,\overline{\St_{\gauta}(\xi)}, \overline{\nu})$, for some measure~$\overline{\nu}$ concentrated on the closure~$\overline{\St_{\gauta}(\xi)}$?  (Problem~8.2 in~\cite{DynSubgroup}). This problem has been studied by Y. Vorobets in the special case of the Grigorchuk group~\cite{vorobets}. He showed that for this group it is possible to reconstruct the action of the general dynamical system on the boundary starting from the study of one orbit. His method uses the study of the map $\SSt$ that associates to any point in the boundary of the tree its stabilizer subgroup in the automaton group. \medskip

Motivated by these ideas, we examine the dynamical and algorithmic implications of the continuity of the map~$\SSt$ in the context of automaton groups, and how some combinatorial properties of the generating automaton reflect into the continuity of this map. In particular, we focus our attention on several classes of Mealy automata: the contracting case (see~\cite{Nekra05}), the reversible case, the bireversible case, and finally the case of automata with a sink-state which is accessible from every state (henceforth denoted by~$\SinkA$). We first show that, in general, the measure of the set of the points in which $\SSt$ is not continuous (henceforth called singular points) is zero. In the bireversible case singular points are exactly points with non-trivial stabilizers. This reproves the well known fact that bireversible automata give rise to essentially free actions on the boundary~\cite[Corollary 2.10]{StVoVo2011}. Driven by these facts and the question raised by Grigorchuk and Savchuk in~\cite{GriSa13} regarding the existence of singular points for the action of a bireversible automaton generating an infinite group, we generalize the previous open problem into the study of examples of automaton groups without singular points. In the case of contracting groups we provide a characterization for such automata in terms of languages recognized by B\"uchi automata that we call \emph{stable automata}. In the examples that we present, we show that the situations may be very different. For instance, in contrast with the Hanoi Towers group case, we show that the Basilica group has no singular points. 
\medskip

In both the class of reversible invertible automata and the class~$\SinkA$, by using the notion of helix graph, we reduce this problem to the existence of certain pairs of words, called {\it commuting pairs}, that is, two words, one on the stateset the other one on the alphabet, that commute with respect to the induced actions.
 Using this fact we prove a series of results. For instance, it turns out that the existence of singular points is always guaranteed for the reversible invertible automata that are not bireversible. This shows that, in the class of reversible invertible automata, the core of the problem of finding examples of group automata without singular points is reduced to the class of bireversible automata. We present some necessary conditions for such examples to exist. For instance, the generated group is necessarily fully positive, that is, it is defined by relators that do not contain negative occurrences of the generators. Furthermore, we prove that if a bireversible automaton generates an infinite non-torsion group, then having all stabilizers in the boundary that are torsion groups (like in the situation of not having singular elements) is equivalent to have in the dual automaton all Schreier graphs in the boundary which are either finite, or acyclic multigraphs (just considering the edges without their inverse).
\medskip

The study of commuting pairs also leads to a connection with periodic tessellations of the discrete plane using Wang tilings, and it is also related to the so-called Gromov's problem (the reader is referred to the paper~\cite{KaPa} for more details). This connection has been pointed out to us by~I.~Bondarenko~\cite{Bonda-private}. Using the helix graph one can easily show that any automaton group has a commuting pair, whence the associated tileset has always a periodic tiling. However this commuting pair may involve a trivial  word. For instance in the class~$\SinkA$ there is always a trivial commuting pair involving the sink-state. This fact leads to the notions of non-elementary commuting pair and reduced tileset of an automaton group. Using a result by~\cite{LeGlo14} we first show that the problem of finding non-elementary commuting pairs is undecidable. Further, the notion of non-elementary commuting pair is strictly related to the existence of periodic singular points in the boundary. From this connection, we start a study of the relationship between non-periodic tessellations of the discrete plane, and algebraic and dynamical properties of the associated automaton group. Indeed, we first provide conditions for the associated reduced tileset to tile the discrete plane. Then, we pinpoint the algebraic and dynamical properties that an automaton group from~$\SinkA$ has to possess so that the associated reduced tileset generates just aperiodic tilings. Finally, we characterize  the existence of aperiodic tilings \textit{\`a~la~} Kari-Papasoglu with some properties of the group generated by an automaton  and its set of singular points.

%-------------------------------------------------------------------------------------------------------------------------------------------------
%-------------------------------------------------------------------------------------------------------------------------------------------------
%-------------------------------------------------------------------------------------------------------------------------------------------------
\section{Preliminaries}

\subsection{Mealy automata}

We first start with some vocabulary on words, then introduce our main tool --- Mealy automata.\\
Let $Q$ be a finite set, as usual, $Q^{n}$, $Q^{\leq n}$, $Q^{<n}$, $Q^{\geq n}$, $Q^{*}$, and $Q^{\omega}$  denote respectively the set of words of length $n$,  of length less than or equal to $n$,  of length less than $n$, of length greater than or equal to $n$, of finite length, and the set of right-infinite words on $Q$ .\\
For two words $u,v\in Q^{*}$ with~$u=vv'$ ($u=v'v$) for some~$v'\in Q^{*}$, we say that $v$ is a \emph{prefix}  (\emph{suffix}),
denoted by~$v\le_{p}u$ (respectively, $v\le_{s}u$). \\
If $\xi=x_1x_2\cdots \in Q^{\omega}$, then \(\xi[n]=x_n\) is the
\(n\)-th letter of \(\xi\), and $\xi[:n]=x_1\cdots x_n$ its
initial prefix  of length $n$. Similarly, for $m\le n$ we
denote by $\xi[m:n]=x_{m}\cdots x_n$ the factor of~$\xi$ of length
$n-m+1$ between the $m$-th and the $n$-th letter of~$\xi$, and by $\xi[m:]=x_{m}x_{m+1}\cdots$  its \emph{tail}. Two infinite sequences $\xi, \eta\in Q^{\omega}$ are said to be \emph{cofinal} (written $\xi\cof\eta$) if there exists an integer $k$ such that $\xi[k:]=\eta[k:]$.\\

By $\tilde{Q}=Q\cup Q\inv$ we denote the \emph{involutive set} where $Q\inv$ is the set of \emph{formal inverses} of~$Q$. 
The operator $\inv\colon Q\rightarrow Q\inv$ sending $q\mapsto q\inv$ is extended to an involution on the free monoid~$\wt{Q}^*$ through
$$
1\inv = 1, \;\; (q\inv)\inv=q, \;\; (uv)\inv=v\inv u\inv\;\;\; (q\in
Q;\;u,v\in\wt{Q}^*).
$$
Let $\sim$ be the congruence on $\wt{Q}^*$ generated by the relation set
$\{(qq\inv,1)\mid q \in \wt{Q} \}$. The quotient $F_Q= \wt{Q}^*/\sim$ is the \emph{free group} on $Q$, and let $\sigma\colon \wt{Q}^* \to F_Q$ be the canonical homomorphism.
The set of all \emph{reduced words} on $\wt{Q}^*$ may be compactly written as
$$
R_Q=\tilde{Q}^*\smallsetminus\bigcup_{q\in\tilde{Q}}\tilde{Q}^*qq\inv\tilde{Q}^*.
$$
For each $u\in\wt{Q}^*$, we denote by $\overline{u}\in R_Q$ the (unique) reduced word $\sim$-equivalent to $u$. With a slight abuse in the notation we often identify the elements of~$F_Q$ with their reduced representatives, \emph{i.e.} $\sigma(u)=\overline{u}$; this clearly extends to subsets $\sigma(L)=\oo{L}$ for~$L\subseteq \wt{Q}^{*}$.\bigskip

A \emph{Mealy automaton} is a tuple $\auta = (\QQ,\XX,\xcdot,\xcirc)$ where $\QQ$ and $\XX$ are finite set respectively called the stateset and the alphabet, and $\xcdot$ , $\xcirc$ are functions from $\QQ \times \XX$ to, respectively, $\QQ$ and $\XX$ called the transition and the production function. This automaton can be seen as a complete, deterministic, letter-to-letter transducer with same input and output alphabet or, following \cite{DaRo14}, as a labelled digraph.\\
The graphical representation  is standard (see Fig.~\ref{fig-lamps-one} for instance) and one displays transitions as follows:
\begin{equation*}
q \xrightarrow{a\mid b} p \ \in \auta \quad \iff \quad  q\xcdot a = p,\quad q\xcirc a = b\:. 
\end{equation*}

It can be seen that the stateset $\QQ$ and the alphabet $\XX$ play a symmetric role, hence we can define a new Mealy automaton: the \emph{dual} of the automaton~$\auta=(\QQ,\XX,\xcdot,\xcirc)$ is the automaton~$\dual\auta=(\XX,\QQ,\xcirc,\xcdot)$ where we have the transition $a\mapright{q|p}b$ whenever $q\mapright{a|b}p$ is a transition in $\auta$ (see Fig.~\ref{fig-lamps-one}). \\ 

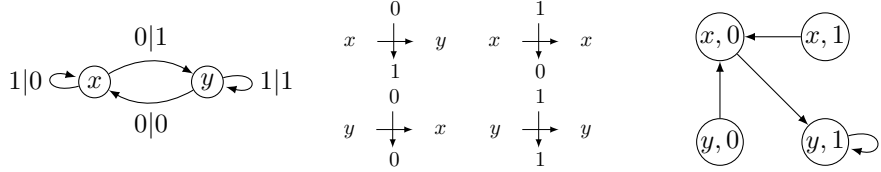
\begin{figure}[ht]%
\centering
\begin{subfigure}[t]{0.33\textwidth}
\begin{tikzpicture}[->,>=latex,node distance=15mm]
	\node[state] (x) {$x$};
	\node[state] (y) [right of=x] {$y$};
	\node (z) [below of=x,node distance=10mm] {};
	\path
		(x)	edge[loop left]		node[left]{\(1|0\)}	(x)
		(x)	edge[bend left]		node[above]{\(0|1\)}			(y)
		(y)	edge[loop right] 	node[right]{\(1|1\)}	(y)
		(y)	edge[bend left]		node[below]{\(0|0\)}	(x);
\end{tikzpicture}
\end{subfigure}
\begin{subfigure}[t]{0.32\textwidth}
\scalebox{.8}{\begin{tabular}[b]{cc}
	{$\begin{array}[b]{ccc}
		& 0            	&    \\
x		& \lacroix    	& y \\
		& 1	 			&            		
\end{array}$} 
&
{$\begin{array}[b]{ccc}
		& 1            	&    \\
x		& \lacroix    	& x \\
		& 0	 			&            		
\end{array}$}\\
	{$\begin{array}[b]{ccc}
		& 0            	&    \\
y		& \lacroix    	& x \\
		& 0	 			&            		
\end{array}$}
 & 
 {$\begin{array}[b]{ccc}
		& 1            	&    \\
y		& \lacroix    	& y \\
		& 1	 			&            		
\end{array}$}
	\\
\end{tabular}
}
\end{subfigure}
\begin{subfigure}[t]{0.32\textwidth}
\begin{tikzpicture}[->,>=latex,node distance=14mm]
	\node[state,inner sep= .3pt] (x0) {$x,0$};
	\node[state,inner sep= .3pt] (y0) [below of=x0] {$y,0$};
	\node[state,inner sep= .3pt] (x1) [right of=x0]{$x,1$};
	\node[state,inner sep= .3pt] (y1) [right of=y0] {$y,1$};
		\node[node distance=10mm] (z0) [left of=x0] {};
		\node[node distance=10mm] (z1) [right of=x1] {};
	\path
		(x0)	edge[]				(y1)
		(y1)	edge[loop right]				(y1)
		(x1)	edge[]					(x0)
		(y0)	edge[]					(x0);
	\end{tikzpicture}
\end{subfigure}
\caption{The Mealy automaton~$\lamp$ generating the lamplighter group, the set of its cross-transitions,
and its helix graph~${\mathcal{H}}_{1,1}(\lamp)$.}
\label{fig-lamps-two}%
\end{figure}

 For each automaton transition $q~\mapright{a|q\xcirc a}~q\xcdot a$, we associate the \emph{cross-transition} depicted in the following way:
\[\begin{array}{ccc}
		& a            	&    \\
q	& \lacroix    	& q\xcdot a, \\
		& q \xcirc	a 			&            		
\end{array}\]
see also Fig.~\ref{fig-lamps-two}.\\
\medskip

If the functions $\left(\XX \to \XX:~ a \mapsto q \xcirc a \right)_{q \in \QQ}$ are permutations the automaton is said to be \emph{invertible}. On the other hand, when the functions $\left(\QQ \to \QQ: q \mapsto q\xcdot a\right)_{a \in \XX}$  are permutations the automaton is called \emph{reversible}. Note that when an automaton is reversible its dual is reversible and the other way around. Mealy automata that are both invertible and reversible are called reversible invertible automata, or RI-automata for short. Other classes of automata will be described in Section~\ref{section:classes}.\\

A Mealy automaton $\auta=(\QQ,\XX,\xcdot,\xcirc)$ defines inductively an action $\QQ\overset{\xcirc}{\curvearrowright} \XX^{*}$ of~$\QQ$ on $\XX^{*}$ by
$$
q\xcirc (a_{1}\cdots a_{n})=(q\xcirc a_{1})\left((q\xcdot a_{1})\xcirc (a_{2}\cdots a_{n})\right)\:,
$$
that can also be depicted by a cross-diagram  by gluing cross-transitions  (see Glasner and Mozes ~\cite{GlaMoz05}, or ~\cite{aklmp12}) representing the action of a word of states on a word of letters (or vice-versa): 
\[\begin{array}{ccccc}
		& a_1            	&            		& a_2 \cdots a_n       			& \\
q		& \lacroix    	& q \xcdot a_1 	& \lacroix				&(q \xcdot a_1) \xcdot (a_2 \cdots a_n)  \\
		& q \xcirc a_1 	&            		&(q \xcdot a_1) \xcirc (a_2 \cdots a_n)  
\end{array} \: .\]

In a dual way, this Mealy automaton defines also an action  $\QQ^{*}\overset{\xcdot}{\curvearrowleft} \XX$. Both actions naturally extend to words, respectively in $\QQ^{*}$ and $\XX^{*}$ with the convention 
\[hg \xcirc a  = h\xcirc (g \xcirc a) \text{ and } g\xcdot ab = (g\xcdot a ) \xcdot b \:.\]

\bigskip

In addition to these descriptions of a Mealy automaton, we are going to use another visualization,
the \emph{helix graphs} (introduced in \cite{aklmp12}).
The \emph{helix graph}~$\mathcal{H}_{n,k}$ of a Mealy automaton $\auta = (\QQ,\XX,\xcdot,\xcirc)$  is the
directed graph with nodes \(\QQ^n\times\XX^k\) and arcs \((u,v)
\longrightarrow \bigl(u\xcdot v, u\xcirc v\bigr)\) for all \((u,v) \in \QQ^n \times \XX^k\) (see Fig.~\ref{fig-lamps-two}).

\subsection{Automaton groups}
From the algebraic point of view, the action
$\QQ^{*}\overset{\xcirc}{\curvearrowright} \XX^{*}$ gives rise to a
semigroup $\semig{\auta}$ generated by the endomorphisms
$q\in \QQ$ of the  regular rooted tree identified with~$\XX^{*}$ defined by
$q: u\mapsto q\xcirc u$ for~$u\in \XX^{*}$.\\
Groups generated by invertible automata play an  important role in group theory (for more
details we refer the reader to~\cite{Nekra05}). In this
framework all the maps~$q: u\mapsto q\xcirc u$, $q\in \QQ$, are automorphisms
of the  regular rooted tree~$\XX^{*}$, and the group generated by these
automorphisms is denoted by $\gauta$ (with identity~$\id$).
Note that the actions $\QQ^{*}\overset{\xcirc}{\curvearrowright} \XX^{*}$
and~$\QQ^{*}\overset{\xcdot}{\curvearrowleft} \XX^{*}$ extend naturally
to the actions $\gauta\overset{\xcirc}{\curvearrowright} \XX^{*}$
and~$\gauta\overset{\xcdot}{\curvearrowleft} \XX^{*}$, respectively.\\
There is a natural way to factorize these actions using the wreath product~\cite{Nekra05,BarSil}.
\medskip

\begin{figure}[ht]%
\centering
\begin{subfigure}[b]{0.312\textwidth}
\scalebox{.9}{
\begin{tikzpicture}[->,>=latex,node distance=14mm]
	\node[state] (x) {$x$};
	\node[state] (y) [right of=x] {$y$};
	
	\node[yshift=-8ex] at ($ (x) !.5! (y) $) (l) {$\lamp$};	
	\path
		(x)	edge[loop left]		node[left]{\(1|0\)}	(x)
		(x)	edge[bend left]		node[above]{\(0|1\)}			(y)
		(y)	edge[loop right] 	node[right]{\(1|1\)}	(y)
		(y)	edge[bend left]		node[below]{\(0|0\)}	(x);
\end{tikzpicture}
}\end{subfigure}
\begin{subfigure}[b]{0.345\textwidth}
\scalebox{.9}{
\begin{tikzpicture}[->,>=latex,node distance=14mm]
	\node[state,inner sep= .3pt] (xx) {$x^{-1}$};
	\node[state,inner sep= .3pt] (yy) [right of=xx] {$y^{-1}$};
		\node[yshift=-8ex] at ($ (xx) !.5! (yy) $) (l) {$\lamp^{-1}$};	
	\path
		(xx)	edge[loop left]		node[left]{\(0|1\)}	(xx)
		(xx)	edge[bend left]		node[above]{\(1|0\)}			(yy)
		(yy)	edge[loop right] 	node[right]{\(1|1\)}	(yy)
		(yy)	edge[bend left]		node[below]{\(0|0\)}	(xx);
	\end{tikzpicture}
}\end{subfigure}
\begin{subfigure}[b]{0.312\textwidth}
\scalebox{.9}{
\begin{tikzpicture}[->,>=latex,node distance=14mm]
	\node[state,inner sep= .3pt] (xx) {$0$};
	\node[state,inner sep= .3pt] (yy) [right of=xx] {$1$};
	\node[yshift=-8ex] at ($ (xx) !.5! (yy) $) (l) {$\dual\lamp$};	
	\path
		(xx)	edge[loop left]		node[left]{\(y|x\)}	(xx)
		(xx)	edge[bend left]		node[above]{\(x|y\)}			(yy)
		(yy)	edge[loop right] 	node[right]{\(y|y\)}	(yy)
		(yy)	edge[bend left]		node[below]{\(x|x\)}	(xx);
	\end{tikzpicture}
}\end{subfigure}
\caption{The lamplighter automaton~$\lamp$, its inverse automaton~$\lamp^{-1}$, and its dual automaton~$\dual\lamp$.}
\label{fig-lamps-one}%
\end{figure}
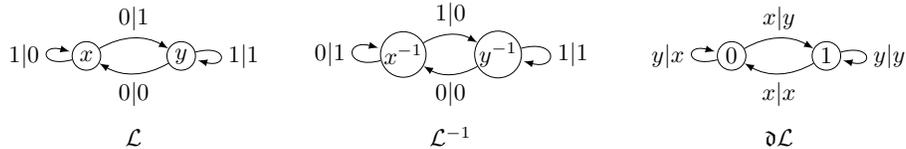

Let $\auta=(\QQ,\XX,\xcdot,\xcirc)$ be an invertible Mealy automaton. The inverse of the automorphism~$q$ is
denoted by~$q^{-1}\in \QQ^{-1}=\{q^{-1}:q\in \QQ\}$. There is an explicit
way to express the actions of the inverses by considering
the inverse automaton~$\auta^{-1}$ having $\QQ^{-1}$ as stateset,
and a transition $q^{-1}\mapright{b|a}p^{-1}$  whenever $q\mapright{a|b}p$ is a transition in~$\auta$ (see Fig.~\ref{fig-lamps-one}).

The action of the group~$\gauta$ on $\XX^{*}$,
in case $\auta$ is invertible (or of the semigroup~$\semig{\auta}$
in a more general case), may be naturally extended on the
\textit{boundary}~$\XX^{\omega}$ of the tree. \\
This action gives rise to the so-called orbital graph. In general, given a finitely generated
semigroup $S$, with set of generators~$Q$, that acts on the left of a
set~$X$ according to  $S\overset{\circ}{\curvearrowright} X$, if $\pi:Q^{*}\rightarrow S$ denotes the canonical map, then
the \emph{orbital graph} $\Gamma(S,Q,X)$ is defined as the $Q$-digraph
with set of vertices~$X$, and there is an edge $x\mapright{a}y$
whenever $\pi(a)\circ x=y$. When we want to pinpoint the connected component
containing the element $y\in X$ we use the shorter notation~$\Gamma(S,Q,X,y)$ instead of~$\left (\Gamma(S,Q,X),y\right)$. Note
that in the realm of groups, this notion corresponds to the notion of
Schreier graph. In particular for an invertible Mealy automaton~$\auta=(\QQ,\XX,\xcdot,\xcirc)$ and a word $v\in \XX^{\ast}\sqcup \XX^{\omega}$, if
\[\St_{\gauta}(v)=\{g\in \gauta
\ : \ g\xcirc v=v\}\] is the stabilizer of~$v$, the Schreier graph
$\Sch(\St_{\gauta}(v), \wt{\QQ})$ corresponds
to the connected component pinpointed by $v$ of the orbital graph:

\[
\Sch(\St_{\gauta}(v), \wt{\QQ})\simeq\Gamma(\gauta,\wt{\QQ},\XX^{*}\sqcup \XX^{\omega},v)
\:.\] 
 This simply corresponds to consider the orbit of~$v$ as the vertex set and the edges given by the action of the generators of the group (in our context the state of the generating automaton). 

Henceforth, when the automaton group is clear from the context we will use the 
more compact notation~$\Sch(v)$ when we deal with~$\Sch(\St_{\gauta}(v), \wt{\QQ})$.

%-------------------------------------------------------------------------------------------------------------------------------------------------
%-------------------------------------------------------------------------------------------------------------------------------------------------
%-------------------------------------------------------------------------------------------------------------------------------------------------
\section{The considered classes}\label{section:classes}

Throughout the paper we focus mainly on four classes of automata: \emph{contracting} automata, \emph{reversible}  automata, \emph{bireversible}  automata, and automata \emph{with a sink}.

\subsection{Contracting automata}
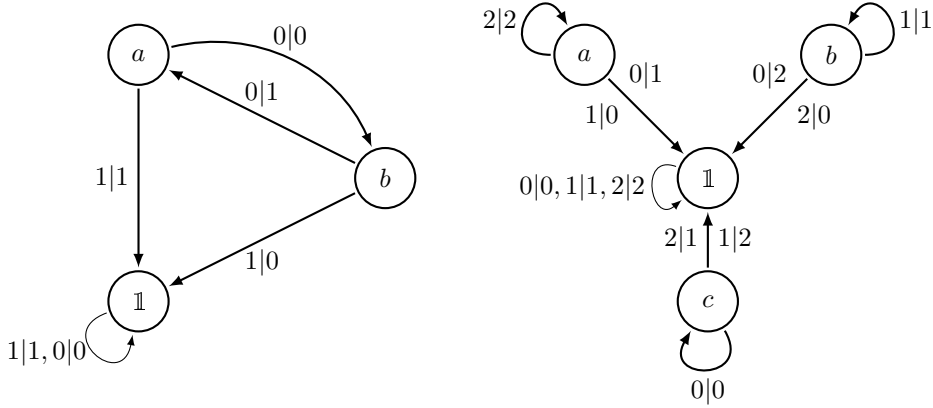
\begin{figure}[h!]
		\begin{center}
			\begin{tikzpicture}[>=latex, shorten >=1pt, shorten <=1pt,scale=.82]
				\tikzstyle{normal_node}= [draw,circle,inner sep=0pt,thick,minimum size=0.8cm]
				\draw (0,0) node [normal_node] (1) {$a$};
				\draw (4,-2) node [normal_node] (0) {$b$};
				\draw (0,-4) node [normal_node] (2) {$\id$};
				\begin{scope}[xshift=35]
					\draw (6,0) node [normal_node] (0') {$a$};
					\draw (10,0) node [normal_node] (1') {$b$};
					\draw (8,-4) node [normal_node] (2') {$c$};		
					\draw (8,-2) node [normal_node] (3') {$\id$};				
				\end{scope}
				\path (0)
				(1) edge[->, thick] node[left] {$1|1$} (2)
				(0) edge[->, thick] node[above] {$0|1$} (1)
				(1) edge[->, bend left=40, thick] node[above] {$0|0$} (0)
				(0) edge[->, thick] node[below] {$1|0$} (2)
				(2) edge[->, out=200, in=260, distance=1cm] node[left] {$1|1, 0|0$} (2)				
				(2') edge[->, thick] node[right] {$1|2$} node[left] {$2|1$} (3')
				(1') edge[->, thick] node[inner sep=10pt, above] {$0|2$} node[inner sep=10pt, right] {$2|0$} (3')
				(0') edge[->, thick] node[inner sep=10pt, above] {$0|1$} node[inner sep=10pt, left] {$1|0$} (3')				
				(0') edge[->, out=180, in=120, distance=1cm, thick] node[left] {$2|2$} (0')
				(2') edge[->, out=-60, in=-120, distance=1cm, thick] node[below] {$0|0$} (2')
				(3') edge[->, out=160, in=220, distance=0.7cm] node[left] {$0|0, 1|1, 2|2$} (3')		
				(1') edge[->, out=0, in=60, distance=1cm, thick] node[right] {$1|1$} (1');			
			\end{tikzpicture}
		\end{center}
		\caption{The contracting automata generating the Basilica group (on the left) and the Hanoi Towers group~$H^{(3)}$ (on the right).} \label{Fig: Bas+Hanoi}
\end{figure}

The notion of contracting automata has been introduced by V. Nekrashevych in~\cite{Nekra05}. For a Mealy automaton $\auta=(\QQ,\XX,\xcdot,\xcirc)$, the group~$\gauta$ is said to be \emph{contracting} if there exists a finite set~$\nucleus\subset\gauta$ such that, for any~$g\in\gauta$ there exists an integer~$n=n(g)$ such that $g\xcdot v\in \nucleus$, for any~$v\in \XX^{\geq n}$. The set~$\nucleus$ is called a \emph{nucleus} of~$\gauta$. Graphically it means that, from any element of the group, a long enough path leads to the nucleus. This enables the construction of a finite automaton~$\auta_{\nucleus}$ with stateset~$\nucleus$,  alphabet~$\XX$, and transitions $g\vvlongmapright{a|g\xcirc a}g\xcdot a$. By extension, an automaton generating a contracting group is said to be \emph{contracting} itself. Examples of such automata are depicted on Fig.~\ref{Fig: Bas+Hanoi} and~\ref{fig: automaton N}.\\
 Note that if the contracting automaton $\auta$ has a sink-state, this state necessarily belongs to~$\nucleus$.\\
The importance of the notion of contracting automata  refers to the beautiful and surprising connection with complex dynamics established by V. Nekrashevych~\cite{Nekra05}.\\
  With every contracting group one may associate a topological space called \emph{limit space}, that is encoded by the set of left-infinite words on $\XX$ modulo the equivalence relation given by the action of the nucleus, \emph{i.e.}, two left infinite sequences $\xi=\cdots\xi_n\xi_{n-1}\cdots\xi_1$ and $\eta=\cdots\eta_n\eta_{n-1}\cdots\eta_1$ are equivalent if for any $n\geq 1$ there exists~$g_n\in \nucleus$ satisfying $g_n\xcirc \xi_n\xi_{n-1}\cdots\xi_1=\eta_n\eta_{n-1}\cdots\eta_1$. It turns out that the iterated monodromy group~$\IMG(f)$ of a post-singularly finite rational function is contracting and its limit space is homeomorphic to the Julia set of~$f$. This discovery puts in strict relation the dynamics of the map $f$ and the algebraic properties of~$\IMG(f)$. As an example, this powerful correspondence has allowed L. Bartholdi and V. Nekrashevych to solve a classical problem in complex dynamics, the so-called Hubbard Twisted Rabbit Problem by using algebraic methods~\cite{BaNe06}.

\medskip

\subsection{(Bi)reversible automata}
The classes of reversible and bireversible Mealy automata are also interesting. We recall that a Mealy automaton is \emph{reversible} whenever each input letter induces a permutation of the stateset, \emph{i.e.} simultaneous transitions $q\mapright{a|b}p$ and $q'\mapright{a|c}p$ are forbidden. In the context of groups, we are especially interested in \emph{reversible invertible} automata (called henceforward $RI$-automata).\\
Moreover such a reversible automaton is \emph{bireversible} if in addition each output letter induces a permutation of the stateset, \emph{i.e.} simultaneous transitions  $q\mapright{a|b}p$ and  $q'\mapright{c|b}p$ are also forbidden. In this case it is necessarily invertible. \\

An interesting feature of an $RI$-automaton is that the dual of such an automaton is still reversible and invertible. \\
The following lemma, which will be useful later, may be easily deduced from~\cite{StVoVo2011} or~\cite[Theorem~2]{DaRo14}.

\begin{lemma}\label{lemma: RI/birev}
Let $\auta$ be an $RI$-automaton with~$\gauta\simeq F_{\QQ}/N$. Then, the following facts hold:
\begin{enumerate}[label=(\roman{enumi})]
\item\label{lemm: RI1} if $g\in \QQ^{\ast}$ is such that $g\xcdot a\in N$ for some~$a\in \XX^*$, then $g\in N$;
\item \label{lemm: RI2} let $g,g',g''\in \QQ^{*}$ with $g\xcdot a=g'$ for some~$a\in \XX^{*}$, then there is a $h\in \QQ^{*}$ satisfying~$(hg)\xcdot a=g''g'$.
\end{enumerate}
Furthermore, if $\auta$ is  bireversible then $\QQ$ can be replaced by $\wt\QQ$ in \ref{lemm: RI1}.
\end{lemma}

\begin{proof}
Let us prove point~\ref{lemm: RI1}. First note that any state reachable from a
state in~\(N\) also belongs to~\(N\). Indeed let \(g'\in N\) and $a',b'\in \XX^*$, we have:
\[
\begin{array}{c c c c c c}

	&	a'		&		& b'	&	\\
g'	&	\lacroix	&	g'\xcdot a'	&\lacroix	& \\
	&	a'		&		& b'	&	\\
\end{array}
\]
So for any word $b' \in \XX^*$, $g'\xcirc a'b'=a'b'$, hence $(g'\xcdot a')\xcirc b'=b'$,
\emph{i.e.}, $g'\xcdot a'\in N$.

 Let $g'=g\xcdot a \in N$. By the reversibility of the automaton
 there is an $ a' \in \XX^*$ such that $ g'\xcdot a'=g$ and
 we can conclude that \(g\) belongs to \(N\).
 
Property~\ref{lemm: RI2} follows by observing that by reversibility there is an $h\in \QQ^{*}$ such that $h\xcdot (g\xcirc a)=g''$ holds (remember that the words in $\QQ^{*}$ are written from right to left). In terms of cross-diagramms we obtain:\\
\[
\begin{array}{c c l}

	&	a			&	\\	
g	&	\lacroix	&	g'=g\xcdot a \\
	&	g\xcirc a	&	\\
h	&\lacroix		&	g'' = h\xcdot (g\xcirc a)\\
	& 			&	
\end{array}
\:.
\]
The last statement follows by applying~\ref{lemm: RI1} to $\auta\sqcup \auta^{-1}$ that is reversible by the bireversibility of~$\auta$.
\end{proof}

The structure of the groups generated by reversible or bireversible automata is far from being understood. For instance, for a long time the only known examples of groups generated by bireversible automata were finite, free, or free products of finite groups~\cite{VoVo07, VoVo10, Nekra05}. Recently examples of bireversible automata generating non-finitely presented groups have been exhibited in~\cite{BoDaRo} and in~\cite{KiPiSa15a,SavSid}. In this regard, we now provide an embedding result of any group generated by a bireversible automaton whose dual does not generate a free group, into the outer automorphism group of a free group of infinite rank. This fact may give some extra insight on the kind of groups that are generated considering bireversibility.

\begin{proposition}\label{prop-embedfree}
Let $\auta$ be a bireversible automaton such that $\gauta$ is infinite and $\gdauta$ is not free. Then, there is a monomorphism
$$
\phi:\gauta\hookrightarrow \Out(F_{\infty}).
$$
\end{proposition}

\begin{proof}
Let $\auta = (\QQ,\XX,\xcdot, \xcirc)$.  By~\cite[Theorem~4]{DaRo15}, if we consider the enriched automaton~$\auta^{-}$ obtained from~$\auta$ by adding the edge~$p\mapright{a^{-1}|b^{-1}}q$ for any edge~$q\mapright{a|b}p$ of~$\auta$, then $\group{\auta^{-}}\simeq\gauta$. By~\cite[Theorem~2]{DaRo14} we may express the group $\gdauta$ as the quotient $F_{\XX}/\oo{N}$ where ${N}$
is the maximal subset invariant under the action $\QQ\overset{\xcirc}{\curvearrowright} \wt{\XX}^{*}$. We may regard ${N}$ as a normal subgroup of the free group~$F_{\XX}$, in particular note that $[F_{\XX}:{N}]<\infty$ if and only if $\gdauta$ is finite, and so $\gauta$ is also finite~\cite{Nekra05,SavVor}. Therefore, $[F_{\XX}:{N}]=\infty$, whence ${N}\simeq F_{\infty}$, since ${N}$ is free by Nielsen's theorem.\\ 
 Let us first prove that there is an embedding $\gauta\hookrightarrow \Aut({N})$. By the stability of~${N}$ under the action $\QQ\overset{\xcirc}{\curvearrowright} \wt{\XX}^{*}$ and the invertibility of~$\auta^{-}$, we have that for any $g\in \gauta$ the map $\psi_{g}: n\mapsto g\xcirc n$ for~$n\in{N}$ is a bijection of~${N}$ that is also a homomorphism since %by (\ref{eq: stability}
$w\xcdot n=w$ holds for any $w\in\wt{\QQ}^{*}$ and~$n\in{N}$. Hence, $\psi_{g}\in \Aut({N})$. Furthermore, the map $\phi:\gauta\rightarrow\Aut({N})$ that sends $g$ to~$\psi_{g}$ is a homomorphism since equality~$\psi_{g'g}=\psi_{g'}\xcirc\psi_{g}$ holds. This map is also injective. Indeed, assume that $\psi_{g_{1}}=\psi_{g_{2}}$ for some~$g_{1}, g_{2}\in \gauta$. Since ${N}$ is normal, then for any $u\in \XX^{*}$ there is a reduced element $n\in {N}$ such that $u$ is a prefix of~$n$ (take a suitable conjugate of a reduced non-trivial element of~${N}$). Hence, $\psi_{g_{1}}=\psi_{g_{2}}$ implies that $g_{1}\xcirc u=g_{2}\xcirc u$ for any $u\in \XX^{*}$, hence~$g_{1}=g_{2}$. Hence we get the claim $\phi:\gauta\hookrightarrow\Aut({N})$. We now show that each automorphism $\psi_{g}$ is not inner. Indeed, assume contrary to our claim, that $\psi_{g}(n)=unu^{-1}$, $n\in{N}$, for some reduced non-empty element $u\in{N}$. Then, for any $n\in{N}$ there is an integer $\ell(n)$ such that $g^{\ell(n)}\xcirc n=n$, whence we have
$$
n=\psi_{g^{\ell(n)}}(n)=\psi_{g}^{\ell(n)}(n)=u^{\ell(n)}nu^{-\ell(n)}
$$
from which we get $nu^{\ell(n)}=u^{\ell(n)}n$. We consider the subgroup
generated by~$n$ and~$u^{\ell(n)}$. By Nielsen's theorem this subgroup
is free, and so both $n$ and $u^{\ell(n)}$ belongs to the same cyclic
subgroup $\la h\ra$ for some~$h\in{N}$. In particular, there
is a non-empty prefix $h'\in\wt{\XX}^{*}$ common to both $u$ and
$n$. However, since $n\in{N}$ is arbitrary and ${N}$
is normal by the same argument above each $u\in\wt{\XX}^{*}$ appears as
a prefix of some non-trivial reduced element $n'\in{N}$, a
contradiction. Hence, $\phi:\gauta\hookrightarrow\Out({N})\simeq
\Out(F_{\infty})$.
\end{proof}
Unfortunately, the condition of having a free group of infinite rank appears to be mandatory in Proposition~\ref{prop-embedfree}. Indeed, the next proposition shows that the embedding of an automaton group generated by a bireversible automaton into the group of length preserving automorphisms of the free group~$F_m$ for some~$1<m<\infty$, characterizes finite groups. In what follows we call an element~$\psi\in \Aut(F_m)$ \textit{length preserving} if given any $w\in F_m$ one has $|w|=|\psi(w)|$. We denote by $\Aut_{\ell\!p}(F_m)$ the subgroup of~$\Aut(F_m)$ formed by the length preserving automorphisms.

\begin{proposition}
Let $\auta$ be a bireversible automaton. There is a monomorphism
$$
\phi: \gauta\hookrightarrow \Aut_{\ell\!p}(F_m)
$$
if and only if $\gauta$ is finite.
\end{proposition}
\begin{proof}
Let $\auta= (\QQ,\XX,\xcdot,\xcirc)$.  Suppose $\gauta$ finite and let ${N}$ be the maximal invariant subset for the action $\QQ\overset{\xcirc}{\curvearrowright} \wt{\XX}^{*}$ (as in Proposition~\ref{prop-embedfree}), then $[F_{\XX}:{N}]<\infty$ and so ${N}\simeq F_m$ for some $m$. We proceed as in the proof of the previous proposition to show that $\phi:\gauta\to \Aut_{\ell\!p}(F_m)$ is a monomorphism. Conversely, let consider an embedding~$\phi: \gauta\hookrightarrow \Aut_{\ell\!p}(F_m)$ for some $m$, and let $R=\{x_1,\ldots, x_m\}$ be minimal set of generators of~$F_m$. For any $g\in \gauta$, let $\psi_g$ be the corresponding automorphism in $\Aut_{\ell\!p}(F_m)$. Since $R$ is finite and the automorphisms preserve the length, the set 
$$
\Omega=\bigcup_{g\in \gauta, i=1,\ldots, m}\psi_g(x_i)
$$
is clearly finite. Further, there is a natural homomorphism of~$\gauta$ into~$Sym(\Omega)$. Let us prove that it is actually a monomorphism. Indeed, let $g\neq g'$ in $\gauta$. Then $\psi_g\neq \psi_{g'}$ holds in $\Aut_{\ell\!p}(F_m)$. Since $R$ is a generating set, then $\psi_g(x_{i_k})\neq \psi_{g'}(x_{i_k})$ for some~$x_{i_{k}}\in R$. We deduce~$\gauta\hookrightarrow Sym(\Omega)$, and so $\gauta$ is finite.
\end{proof}
For similar results that link automaton groups defined by bireversible Mealy automata and the group of automorphisms of a free group,
the reader is referred to~\cite{MaNeSu}.

\subsection{Automata with sink}

In the complement of the class of the $RI$-automata there is
another interesting class that, in some sense, represents the opposite case: the class~$\SinkA$ of all the invertible Mealy automata with a sink-state $e$ which is accessible from every state (the index ``a'' standing for accessible). We recall that a \emph{sink-state} of a Mealy automaton is a special state~$e$ such that $e\xcdot a=e$ and~$e\xcirc a=a$ for any~$a\in A$. Note that in this setting the sink-state is unique.\\
The reason we require that the sink-state is accessible from every state will be clear in Section~\ref{sec: wang tilings}.\\
This class is rather broad and it contains many known classes of Mealy automata like automata with polynomial state activity~\cite{Sidki}. Furthermore, in~\cite[Proposition 6]{DaRo14} it is shown that this class is essentially formed by those automata for which every element $g$ in the generated group has a $g$-regular element in the boundary (for the notion of~$g$-regular element see for instance~\cite{Nekra10}).\\
  Moreover, this class is also included into the broader class of \emph{synchronizing} automata for which some results on automaton groups can be found in~\cite{DaRo13}. The connection with synchronizing automata will also be crucial in Section~\ref{sec: wang tilings} in characterizing automata whose associated set of reduced tiles do not tile the plane. In~\cite{DaRo15}  the problem of finding free groups generated by automata in~$\SinkA$ is tackled. Indeed,until recently, all known free automaton groups were generated by bireversible automata.
This led to the question whether or not it is possible to generate a free group by means of automata with a sink-state. In~\cite{DaRo15}  a series of examples of automata from~$\SinkA$ generating free groups is exhibited. However, in this case the resulting free groups do not act transitively on the corresponding tree, so this leaves open the question of finding a free group generated by an automaton from~$\SinkA$ acting transitively on the rooted tree. This problem is also connected with the interesting combinatorial notion of \emph{fragile word} introduced in~\cite{DaRo15}.

%-------------------------------------------------------------------------------------------------------------------------------------------------
%-------------------------------------------------------------------------------------------------------------------------------------------------
%-------------------------------------------------------------------------------------------------------------------------------------------------
\section{Topological properties of the action on the boundary}\label{sec: dynamics}

In this section we describe some topological properties of the action of an automaton group on the boundary of a rooted tree. In particular, we consider the problem of continuity of the map that associates with any point in the boundary the corresponding Schreier graph. We prove that the set of those points where this function is not continuous has zero measure. Moreover, we provide a characterization of contracting automata whose action on the boundary is continuous everywhere. In the reversible case, we prove that examples of automata generating groups with all continuous points in the boundary are necessarily bireversible, and in this case, this condition may be rephrased in terms of triviality of the stabilizers in the boundary.\\

\subsection{Action on the boundary}

Let $\auta = (\QQ,\XX,\xcdot, \xcirc)$ be an invertible automaton and $Sub(\gauta)$ denote the space of all subgroups of~$\gauta$ and let $\Sch(\gauta,\wt\QQ)$ denote the space of marked Schreier graphs of~$\gauta$ (\emph{i.e.} Schreier graphs in which we have chosen a special vertex, the marked vertex) contained in the space of all marked labeled graphs and put $\partial T=\XX^{\omega}$. Both spaces may be endowed with a natural topology (also induced by an opportune metric). We endow the space~$Sub(\gauta)$ with the Tikhonov topology of the space~$\{0,1\}^{\gauta}$ in such a way that any subgroup $H$ may be identified with its characteristic function. Given a finite subset~$F$ of~$\gauta$ the $F-$neighborhood of a subgroup $H$ contains all subgroups $K$ such that $H\cap F=K\cap F$. Roughly speaking we say that two subgroups~$H$ and~$K$ of~$\gauta$ are close if they share many elements. On the other hand, two marked Schreier graphs $\Sch(\xi)$ and $\Sch(\eta)$ are close when the subgraphs given by the balls of large radius around $\xi$ and $\eta$ are isomorphic, and two points~$\xi$ and~$\eta$ in~$\partial T$ are close if they share a long common prefix. Notice that in our notation, $\Sch(\xi)$ corresponds to the graph $\Gamma(\gauta,\wt{\QQ},\XX^{*}\sqcup \XX^{\omega},\xi)$.
\\

Vorobets studied the map
\begin{align*}
F\colon\partial T	&\longrightarrow	\Sch (\gauta, S)\\
	\xi 			&\longmapsto		\Sch(\St_{\gauta}(\xi), \wt\QQ)
\end{align*}
in the case where $\auta$ is the Grigorchuk automaton~\cite{vorobets}. His results may be summarized as follows: the closure $\overline{F(\partial T)}$ of the image of the boundary of the binary tree into the space of marked labeled Schreier graphs consists of a countable set of points (the one-ended boundary graphs) and another component containing all two-ended Schreier graphs. The Grigorchuk group acts on the compact component given by $\overline{F(\partial T)}$ without these isolated points by shifting the marked vertex of the graph and such action is minimal (every orbit is dense) and uniquely ergodic (there is a unique Borel probability measure on this set that is invariant under the action of the group). 

\subsection{Singular points}
Let $\auta  = (\QQ,\XX,\xcdot,\xcirc)$ invertible. We define the map 
\begin{align*}
\SSt\colon\XX^{\omega}	&\longrightarrow	Sub (\gauta)\\
	\xi 				&\longmapsto		\St_{\gauta}(\xi).
\end{align*}
The \emph{neighborhood stabilizer} $\St_{\gauta}^{0}(\xi)$ of~$\xi$ is the set of
all $g\in\gauta$ that fix the point~$\xi$ together with its neighborhood
(that may depend on $g$). One may check that $\St_{\gauta}^{0}(\xi)$ is a
normal subgroup of~$ \St_{\gauta}(\xi)$.\\
 A point $\xi\in \XX^{\omega}$ is called \emph{singular} if the map $\SSt$ is not continuous at~$\xi$. The
set of singular points is denoted by~$\krit$.\\
 The following lemma clarifies the connection between the continuity of the map~$\SSt$ and the dynamics in the boundary.

\begin{lemma}\cite[Lemma 5.4]{vorobets}
$\SSt$ is continuous at the point $\xi$ if and only if the stabilizer of~$\xi$ under the action coincides with its neighborhood stabilizer, \emph{i.e.}: \[\xi \in \krit \Longleftrightarrow \St_{\gauta}^{0}(\xi) \neq \St_{\gauta}(\xi)\:. \]
\end{lemma}

The following lemma characterizes continuous points in terms of restrictions. 

\begin{lemma}\label{lemmarestric}
Let $\auta$ be an invertible automaton and let $\xi$ be an element in $\XX^{\omega}$.  The following are equivalent.
\begin{enumerate}[label=(\roman{enumi})]
  \item \label{restric1}$\xi$ is not singular;
  \item \label{restric2}$\SSt$ is continuous at $\xi$;
  \item \label{restric3} For any
    $g\in \St_{\gauta}(\xi)$ there exists $n$ such that $g\xcdot \xi[:n]=\id$;

\end{enumerate}
\end{lemma}
\begin{proof}
\ref{restric1}$\Leftrightarrow$\ref{restric2}, \ref{restric3}$\Rightarrow$\ref{restric2},  follow from definition.\\
\ref{restric2}$\Rightarrow$\ref{restric3}. Let us prove $\neg$\ref{restric3}$\Rightarrow$$\neg$\ref{restric2}. Suppose that there exists  $g \in \St_{\gauta}(\xi)$ such that $g \xcdot \xi[:n]\neq \id$ for all $n \geq 0$. Then we can find, for any $n$, a letter $x_n \in \XX$ such that $(g\xcdot \xi[:n])\xcirc x_n = x'_n\neq x_n$. Hence if we put $\zeta_n = \xi[:n]x_n\xi[n+2:]\in \XX^{\omega}$ we get $g \xcirc \zeta_n = g \xcirc \xi[:n]x_n\xi[n+2:]  = \xi[:n]x'_n\xi' \neq \xi$.  Hence,  since we can construct a $\zeta_n$ in any neigbourhood of $\xi$, $\St_{\gauta}^{0}(v^{\omega})\neq \St_{\gauta}(v^{\omega})$, and $\SSt$ is not continuous.
\end{proof}

Moreover we can characterize continuous points by looking only to periodic points
\begin{lemma}\label{lemmaperiodiche}
Let $\auta$ be an invertible automaton. The following are equivalent.
\begin{enumerate}[label=(\roman{enumi})]
\item \label{periodiche1} There is no singular point in $\XX^{\omega}$;
\item \label{periodiche2} There is no singular periodic point in $\XX^{\omega}$.
%\item \label{periodiche2} $\SSt$ is continuous at each $v^{\omega}$, $v\in \XX^{*}$.
\end{enumerate}
\end{lemma}
\begin{proof}
 \ref{periodiche1} $\Rightarrow$ \ref{periodiche2}  is obvious, let us prove the converse, by contraposition.
 Assume that $\SSt$ is not continuous at some~$\xi\in \XX^{\omega}$ and let $g\in \St_{\gauta}(\xi)\smallsetminus \St_{\gauta}^{0}(\xi)$. If there exists $k$ such that $g\xcdot \xi[:k]=\id$ then $g$ stabilizes some neighborhood $U$ of~$\xi$ and  is contained in $\St_{\gauta}^{0}(\xi)$. Hence, for any $n\geq 1$, $g\xcdot \xi[:n]$ is a non-trivial element in $\gauta$. The set $\{\ g\xcdot \xi[:k], \ k\in \mathbb{N}\ \}$ is finite, this implies that there exist $m$ and $n$ such that $n>m>0$ and  $g':=g\xcdot \xi[:m]= g\xcdot \xi[:n]\neq \id$. Therefore $g'\xcdot \xi[m+1:n]=g'$ and $g'\xcirc \xi[m+1:n]=\xi[m+1:n]$. Put $v=\xi[m+1:n]$: $g'\in \St_{\gauta}(v^{\omega})$. In order to prove that $g'\not\in \St_{\gauta}^{0}(\xi)$ we notice that, since $g'$ is not trivial, there exists $w\in \XX^{\ast}$ such that $g'\xcirc w=w'\neq w$.  Consider the sequence $w_k:=v^kwv^{\omega}$ for~$k \geq 0$. Clearly, for any neighborhood $U$ of~$v^{\omega}$ there exists $n$ such that $w_n\in U$. But
$$
g'\xcirc w_n= g'\xcirc v^n (g'\xcdot v^n)\xcirc (wv^{\omega})=v^nw'v'\neq w_n
$$
for some~$v' \in \XX^{\omega}$. Therefore $g'\in \St_{\gauta}(v^{\omega})\neq \St_{\gauta}^{0}(v^{\omega})$.
\end{proof}

In the following theorem we prove that the measure of the set~$\krit$ of singular points is zero. For the sake of completeness we recall that a subset of a topological space $X$ is nowhere dense if its closure has an empty interior. A subset is \emph{meager} in $X$ if it is a union of countably many nowhere dense subsets. A Baire space, as $\XX^{\omega}$ with the usual topology, cannot be given by the countable union of disjoint nowhere dense sets. In general the notion of nowhere dense and meager set do not coincide with the notion of zero-measure. In \cite{vorobets} it is proven that $\krit$ is meager. \\
Given $u\in \wt{\QQ}^{*}$, denote by~$\Fix(u)$ the set consisting  in the vertices $\xi\in \XX^{\omega}$ fixed by the action of~$u$. If $w\in \XX^k$ is an element stabilized by $u$ we write $w \in \Fix_{k}(u)$. By using the ideas developed in  \cite{KaSiSt06} we are able to give the following characterization. 

\begin{theorem}\label{theo: not continuous has zero measure}
For any invertible automaton, the set~$\krit$ of singular points has measure zero.
\end{theorem}
\begin{proof}
The proof is heavily based on the ideas contained in Proposition 4.1 and Theorem~4.2 of~\cite{KaSiSt06}. For $u\in \wt{\QQ}^{*}$, $k\ge 1$, consider the following sets:
$$
\chi(u)=\left\{\xi\in \XX^{\omega}: \xi \in \Fix(u)\mbox{ and }\pi(u\xcdot \xi[:j])\neq \id \quad \forall j\ge 0 \right\}\:,
$$
$$
\chi_{k}(u)=\left\{w\in \XX^{k}: w \in \Fix_{k}(u)\mbox{ and }\pi(u\xcdot \xi[:j])\neq \id \quad 0\le j\le k \right\}\:.
$$
By Lemma~\ref{lemmarestric}, the set of singular points is
$$\krit=\bigcup_{u\in \wt{\QQ}^{*}}\chi(u).$$
Let us prove  $\measure(\chi(u))=0$, for any $u\in \wt{\QQ}^{*}$. Since \[\chi(u)=\bigcap_{k\ge 1}(\chi_{k}(u)\XX^{\omega})\:,\] we get:
$$
\measure(\chi(u))=\lim_{k\to \infty}\measure(\chi_{k}(u)\XX^{\omega})=\lim_{k\to \infty}\frac{|\chi_{k}(u)|}{|\XX|^{k}}\:. 
$$
We now show that this limit is $0$. Indeed, as proved below, there is an integer $p$ such that 
\begin{equation}\label{eq:meas}|\chi_{pk}(u)|\le (|\XX|^{p}-1)^{k}, \text{ for all } k\ge 1 \:.\tag{$E_k$}\end{equation}
Let $H_i=\{v\in \wt{\QQ}^{i}: \pi(v)\neq \id\}$. Since $H_{|u|}$ is finite, there is an integer $p$ such that no element of $H_{|u|}$ induces the identity on $\XX^p$: take it for Equation~\eqref{eq:meas}.\\ 
Use an induction on $k \geq 1$. For $k=1$: $\chi_p(u) \subseteq \XX^p$ and $\Fix_p(u) \neq \XX^p$ by the choice of $p$. Suppose that $|\chi_{p(k-1)}(u)|\le (|\XX|^{p}-1)^{k-1}$. Since $u\xcdot h\in H_{|u|}$ for any~$h\in \chi_{p(k-1)}(u)$,  there is a $v\in \XX^{p}$ that is not fixed by $u\xcdot h$, whence \[|\chi_{p\cdot k}(u)|\le |\chi_{p(k-1)}(u)|(|\XX|^{p}-1)\le (|\XX|^{p}-1)^{k}\:.\]
\end{proof}

\subsection{The contracting case}
The property of being contracting allows us to characterize the set $\krit$ of singular points in terms of a language recognized by an automaton.

Let $\auta=(\QQ,\XX,\xcdot, \xcirc)$ be a Mealy automaton.
We define its \emph{stable automaton} as the automaton~$\Prune{\auta}$ on infinite words where for~$a\in \XX$ we have:
\[ q\mapright{a} p\ \in \Prune{\auta}\quad\Longleftrightarrow\quad
q\mapright{a|a}p\  \in \auta\:.\]

\begin{figure}[h!]
		\begin{center}
			\begin{tikzpicture}[>=latex, shorten >=1pt, shorten <=1pt,scale=.8]
				\tikzstyle{normal_node}= [draw,circle,inner sep=0pt,thick,minimum size=0.85cm]
				\draw (0,-2) node [normal_node] (1) {$b$};
				\draw (0,-4) node [normal_node] (2) {$a$};
				\draw (2,0) node [normal_node] (3) {$ba^{-1}$};
				\draw (4,0) node [normal_node] (4) {$ab^{-1}$};
				\draw (6,-2) node [normal_node] (5) {$a^{-1}$};	
				\draw (6,-4) node [normal_node] (6) {$b^{-1}$};
				\draw (3,-3) node [normal_node] (7) {$\id$};												
				\path 
				(1) edge[->, bend left=40, thick] node[right] {$0|1$} (2)
				(2) edge[->, bend left=40, thick] node[left] {$0|0$} (1)
				(5) edge[->, bend left=40, thick] node[right] {$0|0$} (6)
				(6) edge[->, bend left=40, thick] node[left] {$1|0$} (5)
				(3) edge[->, bend left=40, thick] node[above] {$0|1$} (4)
				(1) edge[->, thick] node[inner sep=8pt, above] {$1|1$} (7)
				(2) edge[->, thick] node[inner sep=8pt, below] {$1|0$} (7)
				(5) edge[->, thick] node[inner sep=8pt, above] {$1|1$} (7)
				(6) edge[->, thick] node[inner sep=8pt, below] {$0|1$} (7)
				(3) edge[->, thick] node[inner sep=8pt, left] {$1|0$} (7)
				(4) edge[->, thick] node[inner sep=8pt, right] {$0|1$} (7)
				(7) edge[->, out=240, in=300, distance=1cm] node[below] {$1|1, 0|0$} (7)
				(4) edge[->, bend left=40, thick] node[below] {$1|0$} (3);						
				\begin{scope}[xshift=-.5]						
				\draw (9,-2) node[state,accepting,inner sep=0pt,minimum size=0.8cm]  (1b) {$b$};
				\draw (9,-4) node[state,accepting,inner sep=0pt,minimum size=0.8cm] (2b) {$a$};
				\draw (10.5,0) node[state,accepting,inner sep=0pt,minimum size=0.8cm] (3b) {$ba^{-1}$};
				\draw (12,0) node [state,accepting,inner sep=0pt,minimum size=0.8cm] (4b) {$ab^{-1}$};
				\draw (13.5,-2) node[state,accepting,inner sep=0pt,minimum size=0.8cm] (5b) {$a^{-1}$};	
				\draw (13.5,-4) node [state,accepting,inner sep=0pt,minimum size=0.8cm] (6b) {$b^{-1}$};
				\draw (11.25,-3) node [normal_node] (7b) {$\id$};		
				\end{scope}
				\path 				
				(2b) edge[->, bend left=40, thick] node[left] {$0$} (1b)
				(5b) edge[->, bend left=40, thick] node[right] {$0$} (6b)
				(2b) edge[->, thick] node[inner sep=8pt, above] {$1$} (7b)
				(5b) edge[->, thick] node[inner sep=8pt, above] {$1$} (7b)
				(7b) edge[->, out=240, in=300, distance=1cm] node[below] {$1, 0$} (7b)	;				
			\end{tikzpicture}
		\end{center}
		\caption{The nucleus automaton  (on the left) and the stable automaton (on the right), associated with the Basilica automaton.} \label{fig: automaton N}
\end{figure}
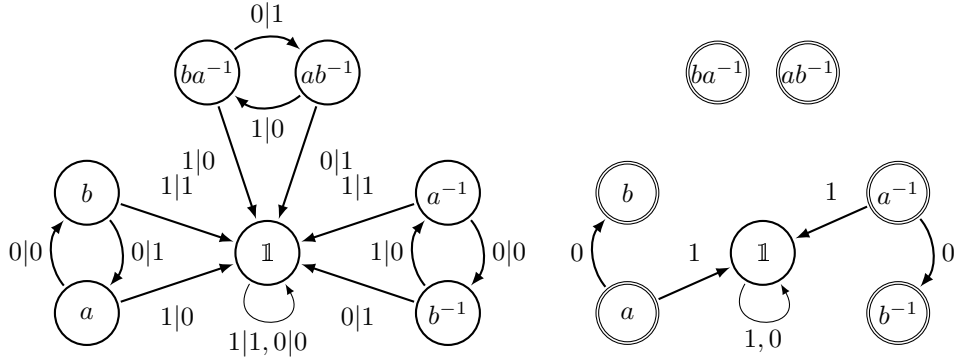

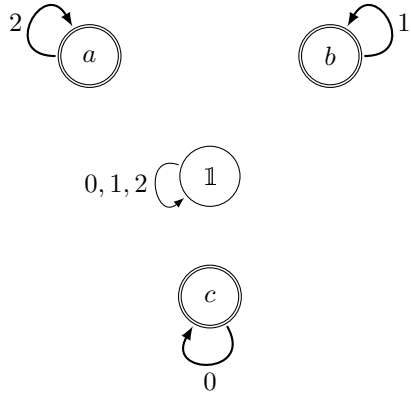
\begin{figure}[h!]
		\begin{center}
			\begin{tikzpicture}[>=latex, shorten >=1pt, shorten <=1pt,scale=.8]
				\tikzstyle{normal_node}= [draw,circle,inner sep=0pt,minimum size=0.8cm]

				\draw (8,0) node[state,accepting,inner sep=0pt,minimum size=0.8cm] (0') {$a$};
				\draw (12,0) node[state,accepting,inner sep=0pt,minimum size=0.8cm] (1') {$b$};
				\draw (10,-4) node[state,accepting,inner sep=0pt,minimum size=0.8cm] (2') {$c$};		
				\draw (10,-2) node [normal_node] (3') {$\id$};

				\path

				(0') edge[->, out=180, in=120, distance=1cm, thick] node[left] {$2$} (0')
				(2') edge[->, out=-60, in=-120, distance=1cm, thick] node[below] {$0$} (2')
				(3') edge[->, out=160, in=220, distance=0.7cm] node[left] {$0, 1, 2$} (3')		
				(1') edge[->, out=0, in=60, distance=1cm, thick] node[right] {$1$} (1');			
			\end{tikzpicture}
		\end{center}
		\caption{The stable automaton~$\Prune{H^{(3)}}$ for the Hanoi Towers automaton.} \label{Fig: stable basilica Hanoi}
\end{figure}

Given a B\"uchi acceptance condition (\emph{i.e.} a set of states that has to be visited infinity often), such an automaton recognizes a language of right-infinite words (see for instance~\cite{Buchi}). See~Figs.~\ref{fig: automaton N} and \ref{Fig: stable basilica Hanoi}.

\begin{theorem}\label{prop:singular_lang}
Let $\auta$ be a invertible automaton admitting a finite nucleus automaton~$\mathcal{N}$. Then the set~$\krit$ of singular points is included in the set of word cofinal with a word in the language recognized by the  B\"uchi automaton $\Prune{\mathcal{N}}$ with every state but the sink-state accepting. Conversely any word recognized by~$\Prune{\mathcal{N}}$ is singular. In particular if the language $\Prune{\mathcal{N}}$ is empty then $\krit = \emptyset$.
\end{theorem}

\begin{proof}
Let $e$ be the sink-state of the nucleus and suppose that there exists
$\xi \in \XX^{\omega}$ such that $\SSt$ is not continuous on
  \(\xi\). From Lemma~\ref{lemmarestric} it follows that there
exists $g \in\St_{\gauta}(\xi)$ such that $g\xcdot\xi[:k] \neq e$ for every $k
> 0$. Since $\gauta$ is contracting there exists $g_n \in \mathcal{N}$ such
that $g_n=g\xcdot \xi[:n] \in \mathcal{N}$. Notice that
$g_{n+k}:=g\xcdot\xi[:n+k] = g_n\xcdot \xi[n+1:n+k] \in
\mathcal{N}\smallsetminus\{e\}$ and $g_{n+k}\xcirc
\xi[n+k] = \xi[n+k]$. Since $\Prune{\mathcal{N}}$ has the same set of states as
$\mathcal{N}$ and admits a transition exactly when $g\xcirc i = i$ holds, there is an infinite path
$$
g_n\mapright{\xi[n]} g_{n+1}\mapright{\xi[n+1]}\cdots g_{n+k}\mapright{\xi[n+k]}\cdots
$$
that avoids $e$. Furthermore since the stateset is finite there is a state that is infinitely visited, hence the run is accepted by $\Prune{\mathcal{N}}$. To prove that $\xi$ is cofinal to a word recognized by $\Prune{\mathcal{N}}$ notice that since $\xi[n:]$ is recognized and $\mathcal{N}$ is finite, there exists an integer~$m$ such that $g\xcdot\xi[:m]=g_{m+1}$ belongs to a strongly connected component~$\mathcal{C}$ of~$\mathcal{N}$. Hence, by reading transition in $\mathcal{C}$ backward,  we can extend $\xi[m:]$ on the left to obtain a word $\xi'=\xi'[1]\cdots \xi'[m-1]\xi[m:]$ and a sequence~$g', g'_1, \ldots, g'_{m}$ of elements of~$\mathcal{C}$   such that there exist~$g'\xcdot \xi'[: k]  = g'_{k+1}$ for~$k \leq m $, $g'\xcdot \xi'[: m] = g_{m+1} $ and  $g'\xcirc \xi'[: k] = \xi  $. Then $\Prune{\mathcal{N}}$ recognizes~$\xi'$, and the points~$\xi'$ and~$\xi$ are cofinal.\\
On the other hand if a run $g_1\mapright{i_1} g_{2}\cdots g_k \mapright{i_k} g_{k+1}\cdots$ is accepted by $\Prune{\mathcal{N}}$ then it is infinite and it does not visit $e$ (since it is a sink, \emph{i.e.} an absorbing state). Hence $g_1$ stabilizes~$\xi = i_1i_2 \cdots$ and satisfies~$\pi(g_1) \neq \id$. So by Lemma~\ref{lemmarestric}, $\SSt$ is not continuous.\end{proof}

Notice that if $\xi \in \XX^{\omega}$ is singular, so is  each element of its orbit. Indeed if $g\xcirc \xi = \xi$ with~$g \neq \id$ then we have $hgh^{-1}\xcirc (h \xcirc \xi) = h\xcirc \xi$ with~ $hgh^{-1}  \neq \id$  for all~$h \in \gauta$. So the previous characterization is exact when the set of words cofinal to a word in the language~$\Prune{\mathcal{N}}$  coincide with the orbit of  words in the language $\Prune{\mathcal{N}}$.\\
In particular this occurs for self-replicating automata. A Mealy  automaton is called \emph{self-replicating} (also called \emph{fractal}) whenever for any word~$u \in \XX^*$ and any~$g \in \gauta$, there exists $h \in \St_{\gauta}(u)$ satisfying~$h\xcdot u = g$. 
\begin{proposition}\label{prop-charac self replic}
Let $\auta$ be a contracting, self-replicating Mealy automaton. Then $\krit$ is the set of words cofinal to a word in the language recognized by $\Prune{\mathcal{N}}$.
\end{proposition}
\begin{proof}
The first inclusion comes directly from Theorem~\ref{prop:singular_lang}. For the other one, let $\xi$ be recognized by $\Prune{\auta}$ and $u\xi[n:]$ be cofinal to $\xi$. 
Since $\xi[n:]$ is recognized by $\Prune{\mathcal{N}}$ it is singular. So there exists~$g \in \gauta$ satisfying~$g\xcirc \xi[n:] = \xi[n:]$ with~$g\xcdot \xi[n: n+k] \neq \id$ for all~$k\geq 0$. Now, since $\auta$ is self-replicating there exists~$h \in \gauta$ satisfying the following cross diagram: 
\[
\begin{array}{c c c c c c}

	&	u		&		& \xi[n:]	&	\\
h	&	\lacroix	&	h\xcdot u = g	&\lacroix	& \\
	&	u		&		& \xi[n:]	&	\\
\end{array}
\]
Hence $u\xi$ is singular.
\end{proof}

Note that, despite of these strong requirements, the class of contracting and self-replicating automata is wide. It contains in particular the Grigorchuk automaton, the Hanoi Tower automaton or the basilica automaton.\\
One can notice that all states in a connected component have same type (accepting or rejecting). In that case the automaton is said to be \emph{weak B\"uchi}~\cite{PerrinPin}. It follows that the language of $\Prune{\mathcal{N}}$ is closed and regular. It is also (infinite) suffix closed.\medskip

Moreover we get informations about the topology of the border $F(\partial T)$. 
\begin{proposition}\label{prop:isol}
Let $\auta$ be a Mealy automaton and $a$ be a letter in its alphabet. If there exists a state $q$ such that $q \xcirc a = a$ and $q\xcdot a = q$,  and $p\xcdot a \neq p$ for~$p \neq q \in \QQ$, then $F(a^{\omega}) = \Sch(a^{\omega})$ is an isolated point in the closure of~$F(\partial T)$. \\
\end{proposition}
\begin{proof}
 Recall that two rooted graphs are close if they have balls of large radius around the roots that are isomorphic. Let us prove~$\St_{\gauta}(a^{\omega})\neq\St_{\gauta}^{0}(a^{\omega})$. The Schreier graph associated with~$a^{\omega}$ contains a loop rooted at~$a^{\omega}$ and labeled by~$q$. Any other element in~$\partial T$ that is not cofinal to~$a^{\omega}$, contains at some position a letter other than~$a$, so that it is not fixed by $q$. This implies that $F(a^\omega)$ is isolated, since no other graph contains a loop labeled by $q$.\\ 
The same argument works for the vertices of the orbit of~$a^{\omega}$, that are the only ones having at finite distance a vertex with a loop labeled by $q$.
\end{proof}
Note that this is still true when one consider words on $\XX^k$ instead of single letters.

We now apply the previous characterizations to two examples of contracting, self-replicating automata, namely the Basilica automaton et the Hanoi Towers automaton.

The Basilica group, introduced in~\cite{GriZu02}, is generated by the automorphisms $a$ and $b$ having the following self-similar form:
$$
a=(b,{\id}), \ \ \ b=(a,{\id})(01),
$$
where $(01)$ denotes the nontrivial permutation of the symmetric group on~$\{0,1,2\}$
and, with a slightly abuse of notation, ${\id}$ denotes the sink-state (see Fig.~\ref{Fig: Bas+Hanoi}). In Fig.~\ref{fig: automaton N} the nucleus automaton associated with this group is presented.

Since the stable automaton accepts no infinite word we obtain:

\begin{corollary}
For the Basilica group, the set~$\krit$ of singular points is empty.
\end{corollary}

We now consider the case of the Hanoi Towers group $H^{(3)}$ (see Fig.~\ref{Fig: Bas+Hanoi}) introduced in~\cite{GriSu06}.
This group is generated by the automorphisms of the ternary rooted tree having the following self-similar form:
$$
a= ({\id},{\id},a)(01) \ \ \ \ b=({\id},b,{\id})(02), \ \ \hbox{and}  \ \ \ c=(c,{\id},{\id})(12),
$$
where $(01)$, $(02)$, and $(12)$ are elements of the symmetric group on~$\{0,1,2\}$. Observe that $a,b,c$ are involutions.

\begin{corollary}\label{theo: hanoi}
For the Hanoi Towers group~$H^{(3)}$, the set~$\krit$ of singular points is a countable set consisting of the (disjoint union of the) orbits of the three points~$0^{\omega}$, $1^{\omega}$, and~$2^{\omega}$. Moreover $F(i^{\omega})=\Sch(i^{\omega})$ is an isolated point in the closure of~$F(\partial T)$ for~$i\in\{0,1,2\}$.
\end{corollary}

{Since $\Sch{(i^{\omega})}$ is isolated in $F(\partial T)$ for~$i\in\{0,1,2\}$, one can ask about the behaviour of the sequence of Schreier graphs of finite words converging to $i^{\omega}$. It turns out that this sequence is converging, but the limit is not a Schreier graph for the Hanoi Towers group.\\ }
Let $\left(\eta_n^i\right)_{n \in \mathbb{N}}$, $i=0,1,2$ be a sequence of elements of $\XX^*$, such that $|\eta_n^i|=n$, converging to $i^{\omega}$. Recall that $\Sch(i^{\omega})$ contains only one loop rooted at~$i^{\omega}$.\\
 Let $\Upsilon_i$, $i=0,1,2$  be the graph obtained from $\Sch(i^{\omega})$ as follows:
\begin{enumerate}
  \item Take two copies of~$\Sch(i^{\omega})$ and let
    $g_i\in\{a,b,c\}$ be the label of the loop at~$i^{\omega}$.
  \item Erase the loop at~$i^{\omega}$ in each copy of~$\Sch(i^{\omega})$.
  \item Join the two copies by an edge labeled $g_i$ and connecting the vertices $i^{\omega}$ of each copy and choose one of these $i^{\omega}$ as marked vertex.
\end{enumerate}

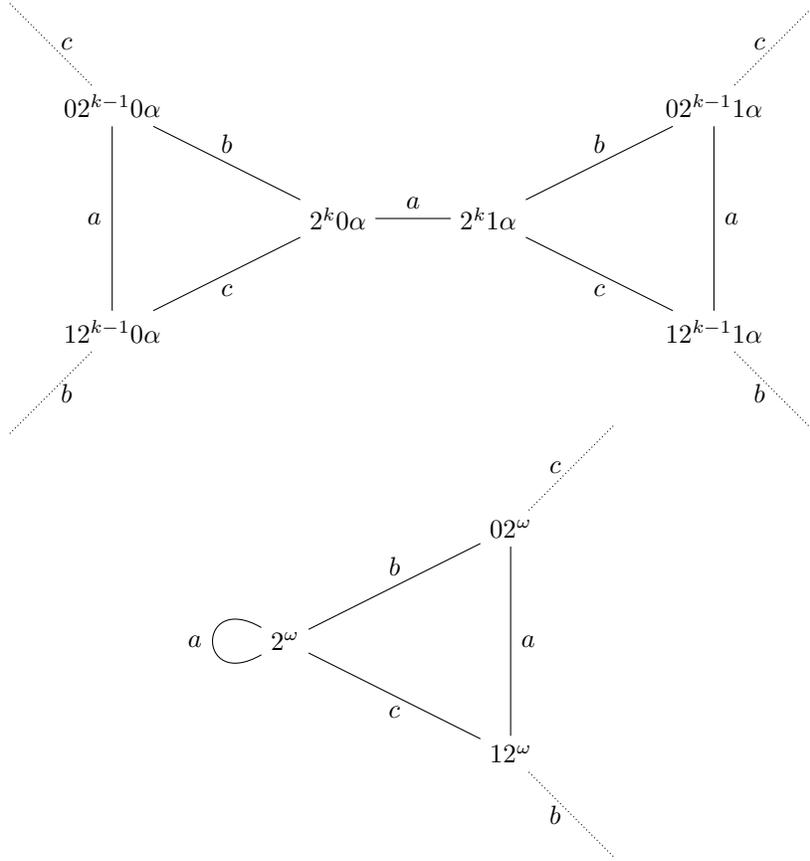
\begin{figure}[h!]
		\begin{center}
			\begin{tikzpicture}[-,>=latex,node distance=15mm]

				\draw  node (20) {$2^k 0\alpha$};
				\draw  node [left of=20] (20') {};				
				\draw  node [left of=20', above of =20'](020) {$02^{k-1}0\alpha$};
				\draw  node [left of=20', below of =20'](120) {$12^{k-1}0\alpha$};
				\draw  node [ left of =020, above of =020,node distance=15mm](020') {};
				\draw  node [left of=120, below of =120,node distance=15mm](120') {};
				
				\draw  node [right of=20,node distance=20mm](21) {$2^k1\alpha$};
				\draw  node [right of=21] (21') {};				
				\draw  node [right of=21', above of =21'](021) {$02^{k-1}1\alpha$};
				\draw  node [right of=21', below of =21'](121) {$12^{k-1}1\alpha$};
				\draw  node [right of=021, above of =021,node distance=15mm](021') {};
				\draw  node [right of=121, below of =121,node distance=15mm](121') {};
				
				\begin{scope}[xshift=-20,yshift=-160]						
				\draw  node (2w) {$2^\omega$};
				\draw  node [right of=2w] (2w') {};				
				\draw  node [right of=2w', above of =2w'](02w) {$02^{\omega}$};
				\draw  node [right of=2w', below of =2w'](12w) {$12^{\omega}$};
				\draw  node [right of=02w, above of =02w,node distance=15mm](02w') {};
				\draw  node [right of=12w, below of =12w,node distance=15mm](12w') {};
				\end{scope}
				
				\path 
				(20) edge[] node[above=0.3pt] {$a$} (21)
				(20) edge[] node[above=0.5pt] {$b$} (020)
				(20) edge[] node[below=0.5pt] {$c$} (120)
				(020) edge[] node[left=0.5pt] {$a$} (120)
				(020) edge[densely dotted] node[right=0.5pt] {$c$} (020')
				(120) edge[densely dotted] node[right=0.5pt] {$b$} (120')

				(21) edge[] node[above=0.5pt] {$b$} (021)
				(21) edge[] node[below=0.5pt] {$c$} (121)
				(021) edge[] node[right=0.5pt] {$a$} (121)
				(021) edge[densely dotted] node[left=0.5pt] {$c$} (021')
				(121) edge[densely dotted] node[left=0.5pt] {$b$} (121')

				(2w) edge[out=210, in=150, distance=1cm] node[left=0.5pt] {$a$} (2w)				
				(2w) edge[] node[above=0.5pt] {$b$} (02w)
				(2w) edge[] node[below=0.5pt] {$c$} (12w)
				(02w) edge[] node[right=0.5pt] {$a$} (12w)
				(02w) edge[densely dotted] node[left=0.5pt] {$c$} (02w')
				(12w) edge[densely dotted] node[left=0.5pt] {$b$} (12w')
				;	
			\end{tikzpicture}
		\end{center}
		\caption{Part of the (finite) Schreier graphs in the Hanoi Towers group~$H^{(3)}$ of $2^k0\alpha$ (above) and of the (infinite) Schreier graph of $2^{\omega}$ (below). See proof of Theorem~\ref{thm-upsilon}.}
		\label{Fig: convergence Schreier graphs}
\end{figure}
\bigskip

\begin{theorem}\label{thm-upsilon}
$F(\eta^i_n)$ converges to $\Upsilon_i$ for $i\in\{0,1,2\}$ as $n\to \infty$.
\end{theorem}
\begin{proof} 
Let $\left( z_n^i\right)_{n \in \mathbb{N}}$ be a sequence of natural numbers such that $z_n^i$ is the position of the
first letter in~$\eta^i_n$ different from $i$. It comes from the
structure of finite Schreier graphs of~$H^{(3)}$ (see Fig.~\ref{Fig: convergence Schreier graphs}) that the balls of
radius~${z_n^i-1}$ in $\Upsilon_i$ rooted at~$i^{\omega}$, and in $\Sch(\eta_n^i)$ rooted at~$\eta_n^i$ are isomorphic. Since $z^i_n$ goes to infinity as $n$ does, we have the assertion.
\end{proof}

It can be shown that the infinite Schreier graphs of the Hanoi Tower group are all one ended, see \cite{BoDaNa16}. Hence, since $\Upsilon_i$ is clearly two ended it follows that it is not an infinite Schreier graph of~$H^{(3)}$. More precisely there is no~$\xi \in \partial T=\XX^{\omega}$ such that the orbital Schreier graph $\Sch(\xi)$ is isomorphic to $\Upsilon_i$, even if the graph is considered non marked.

\begin{figure}[h!]
		\begin{center}
			\begin{tikzpicture}[>=latex, shorten >=1pt, shorten <=1pt,scale=.72]
				\tikzstyle{normal_node}= [draw,circle,inner sep=0pt,thick,minimum size=0.8cm]
				\draw (0,0) node [normal_node] (b) {$b$};
				\draw (2,-2) node [normal_node] (c1) {$c_1$};
				\draw (-2,-2) node [normal_node] (c0) {$c_0$};
				\draw (4,0) node [normal_node] (d1) {$d_1$};
				\draw (-4,0) node [normal_node] (d0) {$d_0$};
				\draw (0,-4) node [normal_node] (a) {$a$};
				\draw (0,-6.5) node [normal_node] (1) {$\id$};
				\path (0)
				(a) edge[->, thick] node[left=-19.5] {$0|1\ \,\,1|0$} (1)
				(b) edge[->, thick] node[below=5,pos=.3] {$1|1$} (c1)
				(c1) edge[->, thick] node[below=5,pos=.7] {$1|1$} (d1)
				(c1) edge[->, thick] node[below=5,pos=.3] {$0|0$} (a)
				(d1) edge[->, thick] node[above] {$1|1$} (b)
				(d1) edge[bend left,->, thick] node[right=4] {$0|0$} (1)
				(b) edge[->, thick] node[below=5,pos=.3] {$0|0$} (c0)
				(c0) edge[->, thick] node[below=5,pos=.7] {$0|0$} (d0)
				(c0) edge[->, thick] node[below=5,pos=.3] {$1|1$} (a)
				(d0) edge[->, thick] node[above] {$0|0$} (b)
				(d0) edge[bend right,->, thick] node[left=4] {$1|1$} (1)
				;		
			\end{tikzpicture}
		\end{center}
		\caption{A twisted version of the Grigorchuk automaton.}\label{Fig: GrigorTwin}
\end{figure}
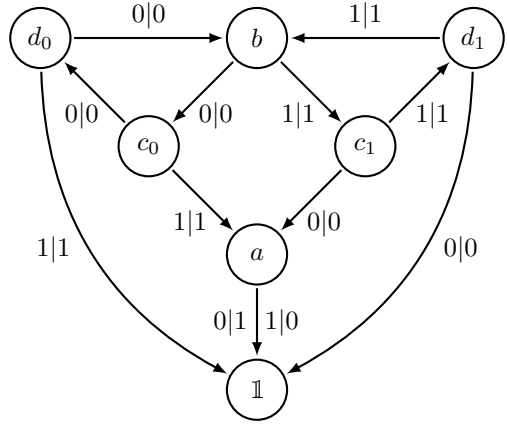\medskip

Note that Grigorchuk automaton can be twisted in order to obtain a contracting automaton generating a fractal group with non countable~$\krit$, see Fig~\ref{Fig: GrigorTwin}. However one can ask if there exist minimal automata where two singular points have isomorphic Shreier graphs, see Problem~\ref{prob:isolated}.
\subsection{The bireversible case}
In this section, we focus on the bireversible case. We start with the following alternative result which shows that, in the class of bireversible automata, the problem of finding examples with all continuous points in the boundary is equivalent to look for automata with all trivial stabilizers in the boundary. 

\begin{proposition}\label{prop: discontinuous point birev}
Let $\auta$ be a bireversible automaton. For any~$\xi\in \XX^{\omega}$, we have: $\xi \in \krit$ $\Leftrightarrow$ $\St_{\gauta}(\xi) \neq \{\id\}$.
\end{proposition}

\begin{proof}
Let $g\in \St_{\gauta}(\xi)$. We have only two possibilities: either there
exists~$n\geq 1$ such that $g\xcdot \xi[:n]=\id$, or for every $n\geq 1$
one has $g\xcdot \xi[:n]\neq \id $. In the first case, since
$\auta$ is bireversible we may apply Lemma~\ref{lemma: RI/birev}, which implies~$g= \id$. On the other hand, suppose $g\xcdot \xi[:n]=g_n\in \gauta\smallsetminus \{\id\}$. We have already remarked that $\SSt : \partial T \longrightarrow Sub (\gauta)$ is not continuous at~$\xi$ if $\St_{\gauta}^{0}(\xi) \neq
\St_{\gauta}(\xi)$. Since $g_{n}\neq \id$, for every $n$ there exists
$w_n\in \XX^{\ast}$ such that $g_n\xcirc w_n\neq w_n$. Let
$\eta^{(n)}:=\xi[:n]w_n \xi[n+|w_{n}|+1:] \in
\XX^{\omega}$. Notice that $\eta^{(n)}\rightarrow \xi$, so that for
every neighborhood~$U$ of~$\xi$ there exists $k=k_U\geq 1$ such that
$\eta^{(k)}\in U$. By hypothesis $g\in \St_{\gauta}(\xi)$ and
$g\xcdot\eta^{(k)}[:k]=g\xcdot \xi[:k]=g_k$. The sequence $\eta^{(n)}$ converges to~$\xi$ since they have the same prefix of length $n$,  so given $U$ we find $k$ such that $\eta^{(k)}\in U$. We have
\begin{eqnarray*}
g\xcirc \eta^{(k)}&=&g\xcirc  \xi[:k]\left((g\xcdot\eta[:k])\xcirc w_k \xi[n+|w_{n}|+1:]\right)\\
&=&\xi[:k](g_k\xcirc w_k)\left((g_k\xcdot w_k)\xcirc \xi[n+|w_{n}|+1:] \right)\\
&\neq&  \eta^{(k)}.
\end{eqnarray*}
This implies that $\St_{\gauta}^{0}(x) \neq
\St_{\gauta}(x)$, which means $x \in \krit$.
\end{proof}
We recall that the action of~$\gauta$ on $\XX^{\omega}$ is \emph{essentially free} if
$$
\measure(\{x\in \XX^{\omega}: \ \St_{\gauta}(x)\neq \{\id\} \})=0,
$$
where $\measure$ is the uniform measure on $\XX^{\omega}$ (see~\cite{DynSubgroup}). One may prove that groups generated by bireversible automata give rise to essentially free actions on $\XX^{\omega}$~\cite{StVoVo2011}. Equivalently, this fact may be deduced by Proposition~\ref{prop: discontinuous point birev} and Lemma~\ref{lemmarestric}. \\
The next proposition characterizes in terms of Schreier graphs those automaton groups having all trivial stabilizers in the boundary,
which by the previous proposition are those with~$\krit=\varnothing$.

\begin{proposition}
Let $\auta$ be a bireversible automaton.
The set~$\krit$ of singular points is empty if and only if any two words in a same orbit cannot be cofinal.
\end{proposition}

\begin{proof}
First, suppose that there exist $\xi\in \XX^{\omega}$ and $\id \neq g\in \gauta$ such that $\eta:=g\xcirc \xi\cof\xi$. Let us prove that there is an element of~$\XX^{\omega}$ whose stabilizer is not trivial. Since $\eta\cof\xi$, then there exists~$N>0$ such that $\xi[N:]=\eta[N:]$. Since $\auta$ is bireversible, from Lemma~\ref{lemma: RI/birev} $g_k:=g\xcdot \xi[:k-1]\neq \id$ for any~$k$. Since $g\xcirc \xi = \eta$, one gets
$$
g_k\xcirc \xi[k:] = (g\xcdot \xi[:k-1]) \xcirc \xi[k:] =  \eta[k:].
$$
In particular, for $k=N$ we get
$$
g_N\xcirc \xi[N:] =g_N\xcirc \eta[N:]= \eta[N:] = \xi[N:].
$$
This implies that $\xi':= \xi[N:]$ is stabilized by $g_N\neq \id$ and so it is a boundary point with a non-trivial stabilizer, so it is a singular point.

On the other hand, let $\gauta$ admit a boundary point $\xi$ with a non-trivial stabilizer. Then, there exists $\id\neq g\in \gauta$ such that $g \xcdot \xi = g$. Since $g \neq \id$ and $\auta$ is bireversible, there exist $ v \in \XX^*$ and $ f \in \QQ^*$, such that $ \: f\xcirc v \neq v$ and $  f \xcdot v = g$. Hence $v\xi $ and $f\xcirc (v\xi)$ are two distinct cofinal words.
\end{proof}

%-------------------------------------------------------------------------------------------------------------------------------------------------
%-------------------------------------------------------------------------------------------------------------------------------------------------
%-------------------------------------------------------------------------------------------------------------------------------------------------
\section{Commuting pairs and dynamics on the boundary}\label{sec: tiling comm}

Theorem~\ref{theo: not continuous has zero measure} from Section~\ref{sec: dynamics} states that the set~$\krit$ of all singular points has measure zero. In this section we are interested in seeking for examples of automata with~$\krit=\varnothing$. Note that by
Proposition~\ref{prop: discontinuous point birev}, in the bireversible
case this property is equivalent to have all trivial stabilizers in
the boundary~\cite{GriSa13}. On the other hand, no example of such
dynamics is known in this class (which seems to be the most difficult case), unless the generated group is finite. In this section we focus our attention on the class of reversible invertible automata.

\subsection{Commuting pairs}
Given an automaton $\auta=(\QQ,\XX,\xcdot,\xcirc)$, we say that $v\in
\XX^{*}$, $u\in \QQ^{*}$ \emph{commute} whenever 
\[
u\xcdot v=u \text{ and } u\xcirc v=v
\ ;\] in this case we say that $(u,v)$ is a \emph{commuting pair}. The previous definition considers words $v\in \XX^{*}, u\in \QQ^{*}$. However, we may consider commuting pairs in~$\auta\sqcup \auta^{-1}=(\wt{\QQ},\XX,\xcdot,\xcirc)$. The importance of commuting pairs stems in the connection with the stabilizers of periodic points on the boundary. For example a classical way to prove that an automaton is not contracting has an interpretation in terms of commuting pairs: if there exists a commuting pair $(u,v)$  such that $u$ has infinite order then the automaton cannot be contracting~\cite{DavEldRee14}. 

We put
\[\St^{+}_{\gauta}(\xi):=\St_{\sgauta}(\xi)
=\pi(\QQ^{*})\cap\St_{\gauta}(\xi)\]
as the set of ``positive'' stabilizers. We have the following proposition that clarifies the connection between commuting pairs and stabilizers.
\begin{proposition}\label{prop: comm stab}
Let $\auta$ be an invertible automaton. If $v\in \XX^{*}, u\in
\wt{\QQ}^{*}$ ($u\in \QQ^{*}$) commute, then $\pi(u)\in
\St_{\gauta}(v^{\omega})$ (respectively, $
\St^{+}_{\gauta}(v^{\omega})$).\\ Conversely, for any $v\in \XX^{*}$, $\pi(u)\in
\wt{\QQ}^{*}$ ($\pi(u)\in \QQ^{*}$) such that $\pi(u)\in \St_{\gauta}(v^{\omega})$
(respectively, $u\in \St^{+}_{\gauta}(v^{\omega})$) then there are integers $n,j\ge 1$ such that $u\xcdot v^{j}, v^{n}$ commute. Moreover if $\auta$ is an RI-automaton then one can take $j=0$, \emph{i.e.}  $u, v^{n}$ commute.

\end{proposition}
\begin{proof}
The first statement follows from $u\xcirc v^{i}=v^{i}$, $u\xcdot v^{i}=u$, for all $i\ge 1$, \emph{i.e.} $\pi(u)\in \St_{\gauta}(v^{\omega})$ ($\pi(u)\in \St^{+}_{\gauta}(v^{\omega})$ in case $u\in \QQ^{*}$). 
On the other hand, since the stateset is finite there exists $k$ such that, $\forall j \geq k$, $u\xcdot v^j$ is in the strongly connected component of $u\xcdot v^k$. Then by strong connectivity there exist~$n$ and~$j \geq k$ such that $(u \xcdot v^j)\xcdot v^n = u \xcdot v^j$. Moreover $(u \xcdot v^j)\xcirc v^n= v^n$, hence  $(u\xcdot v^j, v^n)$ is a commuting pair. In particular if $\auta$ is in addition reversible then every connected component is in fact strongly connected, hence one can choose~$j=0$.
\end{proof}
The helix graph of~$\auta\sqcup \auta^{-1}$ is denoted by $\wt{\mathcal{H}}_{k,n}$. The following proposition gives a way to build all pairs of commuting pairs by looking at the ``labels'' of the cycles of the helix graphs.
\begin{lemma}\label{lem: helix commuting}
Let $\auta$ be an invertible automaton. Let
$$
(u_{0},v_{0})\mapright{}(u_{1},v_{1})\mapright{}\cdots \mapright{}(u_{m},v_{m})\mapright{}(u_{0},v_{0})
$$
be a cycle in the helix graph $\mathcal{H}_{k,n}$ ($\wt{\mathcal{H}}_{k,n}$). Let $u=u_{m}\cdots u_{0}$ and $v=v_{0}\cdots v_{m}$, then $u,v$ commute.\\ Conversely, for any commuting pair $v\in \XX^{*},u\in \QQ^{*}$ ($u\in \wt{\QQ}^{*}$) there is a helix graph $\mathcal{H}_{k,n}$ ($\wt{\mathcal{H}}_{k,n}$) such that $(u,v)\mapright{}(u,v)$.\\ Furthermore, for any $u_{0}\in \QQ^{k}$ ($u_{0}\in\wt{\QQ}^{k}$), $v_{0}\in \XX^{n}$ there is a path 
$$
(u_{0},v_{0})\mapright{}(u_{1},v_{1})\mapright{}\cdots \mapright{}(u_{\ell},v_{\ell})\mapright{}\cdots \mapright{}(u_{\ell},v_{\ell})
$$
in the helix graph $\mathcal{H}_{k,n}$ ($\wt{\mathcal{H}}_{k,n}$), for some~$\ell\ge 0$.
\end{lemma}
\begin{proof}
From the definition of the helix graph we get:
$$
(u_{m}\cdots u_{0})\xcirc (v_{0}\cdots v_{m})= v_{0}\cdots v_{m}, \quad (u_{m}\cdots u_{0})\xcdot (v_{0}\cdots v_{m})= u_{m}\cdots u_{0}
$$
that is, $u,v$ commute. On the other hand, if $u\xcdot v=u$ and $u\xcirc v=v$, then the path~$(u,v)\mapright{}(u,v)$ is a cycle in $\mathcal{H}_{k,n}$ for $k=|u|$, $n=|v|$. The last statement is a consequence of the determinism of an automaton and the fact that any pair~$(u,v)$ has out-degree one, \emph{i.e.}, for~$(u,v)$ there is exactly one pair~$(u',v')$ such that $(u,v)\mapright{}(u',v')$ is an edge in the helix graph.
\end{proof}

\subsection{The reversible-invertible case}
\label{subsec: ri-automata}
In this section we prove that those examples of~$RI$-automata having all
continuous points (for the map $\SSt$) in the boundary (if any) are located in the class of bireversible automata. In other words an $RI$-automaton where the set of singular points is empty is necessarily bireversible. We also show a series of equivalences and connections with the property of having all trivial stabilizers in the boundary.\\

Examples of~$RI$-automata with the aforementioned property are
strictly related to a particular class of groups. Given a group $G$
presented by $\la X|R \ra $, we say that $G$ is \emph{fully positive} 
whenever $R\subseteq X^{*}$, and for any $x\in X$ there is a word~$u\in R$ such that $x$ is a prefix of~$u$. Note that the last property is equivalent to the fact that each $x^{-1}\in X^{-1}$ may be expressed as a positive element $u_{x}\in X^{*}$. These groups have the following alternative (apparently stronger) definition.
\begin{proposition}\label{prop: fully positive}
Let $G$ be a group presented by $\la X|R \ra $, and let $\pi: \wt{X}^{*}\rightarrow G$ be the natural map. The following are equivalent.
\begin{enumerate}[label=(\roman{enumi})]
\item \label{fullpos1}$G$ is fully positive;
\item \label{fullpos2}For any $u\in X^{*}$ there is $v\in X^{*}$ such that $\pi(uv)=\id$ in $G$;
\item \label{fullpos3}For any $u\in X^{*}$ there is $v\in X^{*}$ such that $\pi(vu)=\id$ in $G$.
\end{enumerate}
\end{proposition}

\begin{proof}
\ref{fullpos2}$\Rightarrow$\ref{fullpos1}.  Follows by substituting every  negative occurrence of a generator in a relator by a positive word.  \\
\ref{fullpos1}$\Rightarrow$\ref{fullpos2}. Let us prove the statement by induction on the length~$|u|$. The case~$|u|=1$ follows from the fact that $G$ is fully positive. Therefore, suppose that the statement holds for $|u|<n$ and let us prove it for $|u|=n$. Consider any~$u=au'\in X^{n}$ with~$|u'|=n-1$, for some~$u'\in X^{*}$ and~$a\in X$. By the induction hypothesis there is~$v'\in X^{*}$ such that $\pi(u'v')=\id$. By the definition there is a defining
relation~$ah\in R$ starting with~$a\in X$. Take $v=uv'h\in X^{*}$, then it is easy to check that $\pi(uv'h)=\pi(v)=\id$ holds. Equivalence~\ref{fullpos2}$\Leftrightarrow$\ref{fullpos3} follows by conjugation. 
\end{proof}

Note that all torsion groups are fully positive. The following lemma is a direct consequence of Proposition 15 in~\cite{DaRo14}.

\begin{lemma}\label{lemma: density prop of periodic points}
Let $\auta=(\QQ,\XX,\xcdot,\xcirc)$ be an invertible automaton. Then, for any $\xi\in \XX^{\omega}$
and $\pi(u)\in\St^{+}_{\gauta}(\xi)$ ($\pi(u)\in\St_{\gauta}(\xi)$) there exist $i>j\ge 1$ such that $\pi(u\xcdot
\xi[:j])\in \St^{+}_{\gauta}(\xi[j:i]^{\omega})$ (respectively, $\pi(u\xcdot \xi[:j])\in \St_{\gauta}(\xi[j:i]^{\omega})$).
\end{lemma}

\begin{proof}
Let $\pi(u)\in\St^{+}_{\gauta}(\xi)$, then we have $u\xcirc \xi[:n] = \xi[:n]$ for all $n\ge 1$. Furthermore, by the finiteness of~$\QQ^{|u|}$, and since $\{u\xcdot \xi[:k]\}_{k>0}$ is infinite, there are two indices $i>j\ge 1$ such that $u\xcdot \xi[:i]=u\xcdot \xi[:j]$. From which it follows that
$$
(u\xcdot \xi[:j]) \xcirc (\xi[j:i]^{\omega})=\xi[j:i]^{\omega}
$$
\emph{i.e.}, $\pi(u\xcdot \xi[:j])\in \St^{+}_{\gauta}(\xi[j:i]^{\omega})$. The general case $\pi(u)\in\St_{\gauta}(\xi)$ is treated analogously considering~$\auta\sqcup \auta^{-1}$ instead of~$\auta$.
\end{proof}

From which we derive the two following consequences.

We have the following theorem.

\begin{theorem}\label{theo: alternatives}
Let $\auta=(\QQ,\XX,\xcdot,\xcirc)$ be an $RI$-automaton such that, for all $\xi\in \XX^{\omega}$, $\St^{+}_{\gauta}(\xi)$ contains only torsion elements. Then $\gauta$ is fully positive.
\end{theorem}

\begin{proof}
In particular, $\St^{+}_{\gauta}(v^{\omega})$ contains only torsion elements for all $v\in \XX^{*}$. Take any arbitrary $u_0\in \QQ^{*}$, we show that there is a relation having $u_0$ as a suffix, whence the statement follows by Proposition~\ref{prop: fully positive}. By Lemma~\ref{lem: helix commuting} there is path 
$$
(u_{0},v_{0})\maprightempty\cdots \maprightempty(u_{\ell},v_{\ell})\maprightempty(u_{\ell+1},v_{\ell+1})\maprightempty\cdots \maprightempty(u_{\ell+k},v_{\ell+k})=(u_{\ell},v_{\ell})
$$
in the helix graph $\mathcal{H}_{|u_0|,n}$ for some~$v_0\in \XX^{n}$, $\ell\ge
0$. Since we have a loop around the vertex $(u_{\ell},v_{\ell})$, by Lemma~\ref{lem: helix commuting} $u=u_{\ell+1}\cdots u_{\ell+k}$
and $v=v_{\ell+1}\cdots v_{\ell+k}$ are a commuting pair. Proposition~\ref{prop: comm stab} implies~$\pi(u)\in \St^{+}_{\gauta}(v^{\omega})$, thus $\pi(u^{m})=\id$ for some~$m\ge 1$. Hence $u_{\ell}$ is the suffix of some relation. Since $u_{0}$ and $u_{\ell}$ belong to the same connected component, by using Lemma~\ref{lemma: RI/birev}, we get a relation ending with~$u_0$.
\end{proof}

Note that, in particular, if $\St^{+}_{\gauta}(\xi)=\{\id\}$ for all $\xi\in \XX^{\omega}$, then $\gauta$ is fully positive.

Here are two very  similar results on stabilizers. One when stabilizers are trivial, the other when they contain only torsion elements.

\begin{proposition}\label{prop: positive torsion equivalences1}
\label{prop: density prop of periodic points1}
Let $\auta=(\QQ,\XX,\xcdot,\xcirc)$ be an $RI$-automaton. The following are equivalent.
\begin{enumerate}[label=(\roman{enumi})]
\item \label{posrel1} $\St^{+}_{\gauta}(\xi)=\{\id\}$  for all $\xi\in \XX^{\omega}$;
\item\label{posrel2} $\St^{+}_{\gauta}(v^{\omega})=\{\id\}$  for all $v\in \XX^{\ast}$;
\item \label{posrel3}$\St_{\gauta}(v^{\omega})=\{\id\}$ for all $v\in \XX^{\ast}$;
\item \label{posrel4}$\St_{\gauta}(\xi)=\{\id\}$ for all $\xi\in \XX^{\omega}$.
\end{enumerate}

\end{proposition}
\begin{proof}
\ref{posrel1}$\Leftrightarrow$\ref{posrel2}. The implication is trivial. Let us look at the reciprocal implication: let $\pi(u)\in \St^{+}_{\gauta}(\xi)$. By Lemma~\ref{lemma: density prop of periodic points}, there are integers $i>j\ge 1$ such that $\pi(u\xcdot \xi[:j])\in \St^{+}_{\gauta}(\xi[j:i]^{\omega})$, which is trivial by hypothesis. Hence  $\pi(u\xcdot \xi[:j])=\id$, and so, by reversibility and Lemma~\ref{lemma: RI/birev}, $\pi(u)=\id$.\\
\ref{posrel2}$\Leftrightarrow$\ref{posrel3}. The converse~\ref{posrel3}$\Rightarrow$\ref{posrel2} is trivial. Conversely, if~\ref{posrel2} holds, then by Theorem~\ref{theo: alternatives} we have that $\gauta$ is fully positive. Therefore, for any $u \in \wt{\QQ}$ we can construct a word $u^+ \in \QQ^*$ such that $\pi(u^+)=\pi(u)$. The result follows.\\
\ref{posrel3}$\Leftrightarrow$\ref{posrel4}  follows from~\cite[Proposition 15]{DaRo14}.
\end{proof}

\begin{proposition}\label{prop: positive torsion equivalences2}
\label{prop: density prop of periodic points2}
Let $\auta=(\QQ,\XX,\xcdot,\xcirc)$ be an $RI$-automaton. The following are equivalent.
\begin{enumerate}[label=(\roman{enumi})]
\item \label{posrel1} $\St^{+}_{\gauta}(\xi)$ is formed by torsion elements, for all $\xi\in \XX^{\omega}$;
\item\label{posrel2} $\St^{+}_{\gauta}(v^{\omega})$ is formed by torsion elements for all $v\in \XX^{\ast}$;
\item\label{posrel5} $\St_{\gauta}(v^{\omega})$  is a torsion group for all $v\in \XX^{\ast}$;
\end{enumerate}
Furthermore, if $\auta$ is bireversible, following is also equivalent:
\begin{enumerate}[label=(\roman{enumi}),resume]

\item\label{posrel6} $\St_{\gauta}(\xi)$ is a torsion group for all $\xi\in \XX^{\omega}$.
\end{enumerate}
\end{proposition}
\begin{proof}
\ref{posrel1}$\Leftrightarrow$\ref{posrel2}. The implication is trivial. Let us look at the reciprocal implication: let $\pi(u)\in \St^{+}_{\gauta}(\xi)$. By Lemma~\ref{lemma: density prop of periodic points}, there are integers $i>j\ge 1$ such that $\pi(u\xcdot \xi[:j])\in \St^{+}_{\gauta}(\xi[j:i]^{\omega})$, which contains only torsion elements by hypothesis. Hence,  there is an integer $\ell\ge 1$ such that
$\pi\left((u\xcdot \xi[:j])^{\ell}\right)=\id$. Further, since $\pi(u)\in
\St^{+}_{\gauta}(\xi)$ the following equality:
$$
(u\xcdot \xi[:j])^{\ell}=u^{\ell}\xcdot \xi[:j]
$$
holds for any $j\ge 1$. Thus, by Lemma~\ref{lemma: RI/birev},
we deduce~$\pi(u)=\id$ and we get $\pi(u^{\ell})=\pi\left((u\xcdot \xi[:j])^{\ell}\right)=\id$, \emph{i.e.}, $\pi(u)$ is torsion.
\\
\ref{posrel2}$\Leftrightarrow$\ref{posrel3}. The proof is similar as the one in Proposition \ref{prop: density prop of periodic points1}.\\
In the bireversible case, the  implication \ref{posrel6}$\Rightarrow$\ref{posrel5} is trivial, while the converse 
is proven similarly to  Proposition~\ref{prop: density prop of periodic points1}. Take $\pi(u)\in \St_{\gauta}(\xi)$. By Lemma~\ref{lemma: density prop of periodic points}, there are integers~$i>j\ge 1$ such that $\pi(u\xcdot \xi[:j])\in \St_{\gauta}(\xi[j:i]^{\omega})$. Since the group $\St_{\gauta}(\xi[j:i]^{\omega})$ is torsion, then there is an integer $\ell\ge 1$ such that $\pi\left((u\xcdot \xi[:j])^{\ell}\right)=\id$. Hence, by the bireversibility and Lemma~\ref{lemma: RI/birev}, we get $
\pi(u^{\ell})=\pi(u^{\ell}\xcdot  \xi[:j])=\pi\left((u\xcdot \xi[:j])^{\ell}\right)=\id$.
\end{proof}

In particular, note that the previous propositions imply the existence of a non-trivial ``positive'' stabilizer whenever the action on the boundary for an $RI$-automaton has at least one non-trivial stabilizer.
Let $\auta=(\QQ,\XX,\xcdot,\xcirc)$ be an invertible automaton, and let $\gauta=F_{\QQ}/N$. The set of ``positive relations'' of~$\gauta$ is
$$
\mathcal{P}(\auta)=\QQ^{+}\cap \psi^{-1}(N)
$$
where $\psi:\wt{\QQ}^{*}\rightarrow F_{\QQ}$  is the canonical homomorphism. Note that $\mathcal{P}(\auta)=\varnothing$ implies that $\semig{\auta}$ is torsion-free and therefore infinite. 

\begin{lemma}\label{cor: empty positive relations implies stab}
If $\mathcal{P}(\auta)=\varnothing$, then $\St^{+}_{\gauta}(v^{\omega})\neq \{\id\}$ for some~$v\in \XX^{*}$.
\end{lemma}
\begin{proof}
Let $$
(u_{0},v_{0})\maprightempty\cdots \maprightempty(u_{\ell},v_{\ell})\maprightempty(u_{0},v_{0})
$$ be a cycle in the helix graph $\mathcal{H}_{1,1}$. Then by Lemma~\ref{lem: helix commuting} $(u_0\hdots u_\ell, v_0 \hdots v_\ell)$ is a commuting pair. Hence $\pi(u_0\hdots u_\ell) \in \St^+_{\gauta}((v_0\hdots v_\ell)^\omega)$, and since $\mathcal{P}(\auta) = \emptyset $, we conclude~$\St^+_{\gauta}((v_0\cdots v_\ell)^\omega) \neq \{\id \}$.
\end{proof}

The next theorem shows that either~$\gauta$ or~$\gdauta$ may have all trivial stabilizers in the boundary. We first recall the following proposition.

\begin{proposition}\cite[Corollary 5]{DaRo15}\label{prop: no positive dual}
Let $\auta=(\QQ,\XX,\xcdot, \xcirc)$ be an $RI$-automaton with~$\gauta$ infinite. Then the index
$\left[\gauta:\St_{\gauta}(y^{\omega})\right]$ is infinite for all $y\in \XX^{*}$,
if and only if
$\mathcal{P}(\dual\auta)=\varnothing$.
\end{proposition}

The following simple proposition is technical:

\begin{proposition}\label{prop: alternative with dual}
Let $\auta=(\QQ,\XX,\xcdot, \xcirc)$ be an $RI$-automaton with $\gauta$ infinite. If \/ $\St^{+}_{\gauta}(y^{\omega})= \{\id\}$ for all $y\in \XX^{*}$, we have $\mathcal{P}(\dual\auta)=\varnothing$.
\end{proposition}

\begin{proof}
Suppose $\St^{+}_{\gauta}(y^{\omega})= \{\id\}$ for all~$y\in \XX^{*}$. Let us prove
\[\left[\gauta:\St_{\gauta}(y^{\omega})\right]=\infty\]
for all~$y\in \XX^{*}$. Indeed, if $\left[\gauta:\St_{\gauta}(y^{\omega})\right]<\infty$ for some~$y\in \XX^{*}$, then the following property holds:
\begin{equation}\label{eq: power}
\exists k\le \left[\gauta:\St_{\gauta}(y^{\omega})\right]\mbox{ such that } \forall g\in \St^+_{\gauta}(y^{\omega}):g^{k}=\id.
\end{equation}
In particular, by Proposition~\ref{prop: fully positive} we get that $\gauta$ is fully positive. We claim that every element of~$\gauta$ is torsion. Indeed, since $\gauta$ is fully positive we have that for any $q\in \QQ$ there is a $u_{q}\in \QQ^{*}$ such that $\pi(q^{-1})=\pi(u_{q})$. Now take any $w\in \wt{\QQ}^{*}$, by substituting each negative occurrence $q^{-1}\in \QQ^{-1}$ appearing in~$w$ with~$u_{q}$, we obtain a ``positive'' word $\widehat{w}\in \QQ^{*}$ such that $\pi(\widehat{w})=\pi(w)$. Thus, by (\ref{eq: power}) we have:
$$
\pi(w)^{k}=\pi(\widehat{w})^{k}=\pi(\widehat{w}^{k})=\id,
$$
that is, $\gauta$ is torsion. Being a residually finite group with uniformly bounded torsion, we deduce from~\cite{Zelma90, Zelma91} that $\gauta$ is finite, a contradiction. Hence, we obtain~$\left[\gauta:\St_{\gauta}(y^{\omega})\right]=\infty$ for all $y\in \XX^{*}$, and so, by Proposition~\ref{prop: no positive dual}, $\mathcal{P}(\dual\auta)=\varnothing$. In particular, by Corollary~\ref{cor: empty positive relations implies stab}, $\mathcal{P}(\dual\auta)=\varnothing$ implies $\St^{+}_{D}(z^{\omega})\neq \{\id\}$ for some~$z\in \QQ^{*}$, hence the last statement holds.
\end{proof}

The latter admits a partial converse in the case where there exists an aperiodic element in $\gauta$:
\begin{proposition}
Let $\auta$ be a RI Mealy automaton. If $\gauta$ is not torsion and $\mathcal{P}{(\dual{\auta})}\neq\varnothing$  then 
\[\exists y \in \XX^*, \: \St^+_{\gauta}(y)\neq \{\id\}\]
\end{proposition}

\begin{proof}
Let $u \in \wt\QQ^*$ such that $\pi(u)$ is aperiodic. Since $\mathcal{P}{(\dual{\auta})}$ is not empty there exists $y$ such that $\sigma(y)= \id$. Then in the helix graph $\wt{\mathcal{H}}(\auta)_{|u|,|v|}$, there is a path 
$$
(u,y)\maprightempty(u,y_{1})\maprightempty\cdots \maprightempty(u,y_{k})\maprightempty(u,y_{k+1})\maprightempty\cdots \maprightempty(u,y_\ell) \maprightempty(u,y_k)
\:.$$
Hence $u^{l-k}$ and $y_k\cdots y_l$ commutes. Then $\pi(u^{l-k}) \in \St_{\gauta}({(y_k\cdots y_l)^{\omega}})$ and $\pi(u^{l-k}) \neq \id$. We conclude using Prop.~\ref{prop: positive torsion equivalences1}.
\end{proof}

\begin{theorem}\label{thm: alternative with dual}
Let $\auta=(\QQ,\XX,\xcdot, \xcirc)$ be an $RI$-automaton with $\gauta$ infinite. We have either $\St^{+}_{\gauta}(y^{\omega})\neq\{\id\}$ for some~$y\in \XX^{*}$, or $\St^{+}_{\gaut{\dual{\auta}}}(z^{\omega})\neq\{\id\}$ for some~$z\in \QQ^{*}$.
\end{theorem}

{ Now, we focus on the possible consequences of being bireversible or not.} We first recall the following proposition.

\begin{proposition}\cite
{DaRo15,GoKiPi14}\label{prop: not containing bireversible}
Let $\auta$ be an $RI$-automaton. If $\auta$ does not contain a bireversible connected component, then $\mathcal{P}(\auta)=\varnothing$.
\end{proposition}

From which we have the following proposition.

\begin{proposition}
Let $\auta=(\QQ,\XX,\xcdot, \xcirc)$ be a non bireversible $RI$-automaton. Then $\St^{+}_{\gauta}(v^{\omega})\neq\{\id\}$ for some~$v\in \XX^{*}$. Furthermore, the set~$\krit$ of singular points is not empty.
\end{proposition}

\begin{proof}
The first statement follows by applying  Proposition~\ref{prop: not containing bireversible} and Corollary~\ref{cor: empty positive relations implies stab} to a non bireversible component of $\auta$. Let us prove the last claim of the proposition. We show that $\SSt$ is actually not continuous at~$v^{\omega}$. Indeed, let $\pi(u)\in\St^{+}_{\gauta}(v^{\omega})$ with~$\pi(u)\neq \id$, for some~$u\in \QQ^{*}$. If $\SSt$ were continuous at~$v^{\omega}$ then, by Lemma~\ref{lemmaperiodiche}, $\pi(u\xcdot v^{j})=\id$ for some~$j\ge 1$. Thus, by Lemma~\ref{lemma: RI/birev}, we get $\pi(u)=\id$, a contradiction.
\end{proof}

As an immediate consequence of the previous proposition we obtain the following result.

\begin{corollary}\label{cor: ri-automaton continuous point are bireversible}
If there exists an $RI$-automaton $\auta$ generating a group without singular points, then necessarily $\auta$ is bireversible.
\end{corollary}

Let us have a look at the case of bireversible automata without singular points.

\begin{proposition}\label{prop: bir+all trivial stabilizers properties}
Let $\auta=(\QQ,\XX,\xcdot, \xcirc)$ be a bireversible automaton with $\gauta$ infinite and no singular points. Then for any $u\in \QQ^{*}$ the following are equivalent:
\begin{enumerate}[label=(\roman{enumi})]
\item \label{bir-triv-stab1}  the Schreier graph centered at~$u^{\omega}$ is finite;
\item \label{bir-triv-stab2}  there is an integer $\ell>0$ such that $\pi(u^{\ell})=\id$;
\item \label{bir-triv-stab3}  $\St_{\gdauta}^{+}(u^{\omega})\neq\{\id\}$.
\end{enumerate}
\end{proposition}
\begin{proof}
Recall that, for bireversible automata~$\auta$, having no singular points is equivalent to $\St^{+}_{\gauta}(v^{\omega})= \{\id\}$ for all $v\in \XX^{*}$, by Propositions~\ref{prop: positive torsion equivalences2} and~\ref{prop: discontinuous point birev}.\\
\ref{bir-triv-stab1}$\Leftrightarrow$\ref{bir-triv-stab2} follows from~\cite[Theorem~6]{DaRo15}.\\
\ref{bir-triv-stab3}$\Rightarrow$\ref{bir-triv-stab2}:  Let $v\in \XX^{*}$ such that
  $\sigma(v)\in \St_{\gdauta}^{+}(v^{\omega})$, where $\sigma: \wt{\XX}^{*}\rightarrow \gdauta$ is the natural map. By Proposition~\ref{prop: comm stab} there is an integer $\ell>0$ such that $v$ and $u^{\ell}$ commute, hence $\pi(u^{\ell})\in \St^{+}_{\gauta}(v^{\omega})=\{\id\}$, \emph{i.e.}, $\pi(u^{\ell})=\id$.\\ 
\ref{bir-triv-stab1}$\Rightarrow$\ref{bir-triv-stab3}: If  the Schreier graph centered at~$u^{\omega}$ is finite, then, for any $v\in \XX^{*}$ we get $\sigma(v^{k})\in \St_{\gdauta}(u^{\omega})$ for~$k=[\gdauta: \St_{\gdauta}(u^{\omega})]$. Furthermore, $\sigma(v^{k})\neq \id$ since, by Proposition~\ref{prop: alternative with dual}, $\mathcal{P}(\dual\auta)=\varnothing$, whence $\St_{\gdauta}^{+}(u^{\omega})\neq \{\id\}$.
\end{proof}

The equivalence \ref{bir-triv-stab2}$\Longleftrightarrow$\ref{bir-triv-stab3} of Prop.~\ref{prop: bir+all trivial stabilizers properties} links being torsion and having non-trivial stabilizers. We name it for future references in this paper:
\begin{equation}\label{eq: relations on non-trivial}
\forall u\in \XX^{*}\mbox{ we have }\St_{\gauta}^{+}(u^{\omega})\neq \{\id\}\Longleftrightarrow \sigma(u^{\ell})=\id\mbox{, for some }\ell\ge 1 \tag{TS}
\:.\end{equation}

We recall that an \emph{acyclic multidigraph} is a multidigraph without cycles. We have the following geometrical description in terms of some algebraic conditions.
\begin{proposition}\label{prop: geo description}
Let $\auta$ be an $RI$-automaton. If $\mathcal{P}(\auta)=\varnothing$ and
\eqref{eq: relations on non-trivial}
then the orbital graphs $\Gamma(\gauta,\QQ, \XX^{\omega},v^{\omega})$ of periodic points $v^{\omega}$ for~$v\in \XX^{*}$, are either finite or acyclic multidigraphs.
\end{proposition}
\begin{proof}
Suppose that both $\mathcal{P}(\auta)=\varnothing$ and Condition~(\ref{eq: relations on non-trivial}) hold. Let $u\in \XX^{*}$. By~\cite[Lemma 2]{DaRo15} it follows that every vertex of~$\Gamma(\gauta,\QQ,
\XX^{\omega},u^{\omega})$ is a periodic point $h^{\omega}$ for some~$h\in \XX^{*}$. Assume that there exists $h \in \XX^{*}$ such that $\St_{\gauta}^{+}(h^{\omega})\neq\{\id\}$, then by condition (\ref{eq: relations on non-trivial}) and~\cite[Theorem~6]{DaRo15}, $\Gamma(\gauta,\QQ,
\XX^{\omega},h^{\omega})$ is finite and so is~$\Gamma(\gauta,\QQ,
\XX^{\omega},u^{\omega})$.\\
 Otherwise, $\St_{\gauta}^{+}(p)=\{\id\}$ for each vertex $p$ of~$\Gamma(\gauta,\QQ, \XX^{\omega},u^{\omega})$. Thus, by  condition $\mathcal{P}(\auta)=\varnothing$, we deduce that, for any vertex $p$ of~$\Gamma(\gauta,\QQ, \XX^{\omega},u^{\omega})$ there is no cycle $p\mapright{v} p$ for any~$v\in \QQ^{*}$, \emph{i.e.}, $\Gamma(\gauta,\QQ, \XX^{\omega},u^{\omega})$ is an acyclic multidigraph.
\end{proof}
Note that the converse  of the previous proposition holds if one assumes the existence of an infinite orbital graph rooted at some periodic point. Gathering Theorem~\ref{thm: alternative with dual} and the equivalence~\eqref{eq: relations on non-trivial}, we obtain the following corollary, and so by Proposition~\ref{prop: geo description} a description of the Schreier graphs of the dual in the case $\gauta$ has no singular points.
\begin{corollary}\label{cor: charact}
Let $\auta$ be a bireversible automaton with $\group{\dual\auta}$ infinite and no singular points in the dual: $\forall \xi \in \QQ^\omega,\: \St_{\gdauta}(\xi)= \{ \id\}$. Then \eqref{eq: relations on non-trivial} holds.
\end{corollary}

We can now obtain a lower bound on the growth of the Schreier graphs pointed at periodic points.
\begin{proposition}
Let $\auta$ be a bireversible automaton with $\gauta$ infinite and no singular points. Then for any $u\in \QQ^{*}$ with $\pi(u)$ aperiodic, we have
$$
\forall m\ge 1, \quad
\left[\gdauta:\St_{\gdauta}(u^{m})\right ]> \log_{|\XX|}(m).
$$
\end{proposition}
\begin{proof}
Fix some \(u\in \QQ^*\) with $\pi(u)$ aperiodic. We claim that for any $v\in \XX^{*}$, if
$\sigma(v)\in\St_{\gdauta}^{+}(u^{m})\smallsetminus\{\id\}$, then $m<|\XX|^{|v|}$. For~$0\le i<m$, we have $\sigma(u^{i}\xcirc v)\in \St^+_{\gdauta}(u)$. By the reversibility of~$\auta$, there is an integer $k\leq|\XX|^{|v|}$ such that $u^{k}\xcirc v=v$.
If~$m\geq|\XX|^{|v|}$, then~$k\leq m$ and $\sigma(v)\in \St^{+}_{\gdauta}(u^{\omega})$, contradicting~Corollary~\ref{cor: charact}. Therefore, $m<|\XX|^{|v|}$. In particular, taking any $a\in \XX$, we get that there is an integer $j\le [\gdauta:\St_{\gdauta}(u^{m}) ]$ such that $\sigma(a^{j})\in \St_{\gdauta}(u^{m})$, from which we obtain
$$
\left[\gdauta:\St_{\gdauta}(u^{m})\right ]\ge j> \log_{|\XX|}(m)
$$
and this concludes the proof.
\end{proof}

We obtain the following geometrical characterization.
\begin{theorem}\label{theo: characterization}
Let $\auta$ be a bireversible automaton with
$\gauta$ infinite, with at least one aperiodic element. The following are equivalent.

\begin{enumerate}[label=(\roman{enumi})]
\item \label{characterizationi} 
 $\St_{\gauta}(\xi)$ is a torsion group, for all~$\xi\in \XX^{\omega}$;
 \item \label{characterizationii}
 $\mathcal{P}(\dual\auta)=\varnothing$ and~(\ref{eq: relations on non-trivial}) hold.
 \item \label{characterizationiii}
 the orbital graphs $\Gamma(\gdauta,\XX, \QQ^{\omega},v^{\omega})$ of periodic points are either finite or acyclic multidigraphs.
\end{enumerate}

\end{theorem}
\begin{proof}
\ref{characterizationi}$\Rightarrow$\ref{characterizationii}:
From \ref{characterizationi}, the subgroup~$\St_{\gauta}(y^{\omega})$ is torsion for any~$y\in \XX^{*}$. Let us prove $\mathcal{P}(\dual\auta)=\varnothing$. Assume $\left[\gauta:\St_{\gauta}(y^{\omega})\right]<\infty$ for some $y\in \XX^{*}$. Then, the subgroup
$$
N=\bigcap_{g\in \gauta}g\St_{\gauta}(y^{\omega}) g^{-1}
$$
is a finite index torsion normal subgroup of~$\gauta$. Thus, since $\gauta/N$ is
finite, the group~$\gauta$ is torsion, a contradiction. We deduce
$\left[\gauta:\St_{\gauta}(y^{\omega})\right]=\infty$ for all $y\in \XX^{*}$,
whence by Proposition~\ref{prop: no positive dual},
$\mathcal{P}(\dual\auta)=\varnothing$. Now, to prove (\ref{eq: relations on non-trivial}) we essentially repeat the proof of the equivalence~\ref{bir-triv-stab3}$\Leftrightarrow$\ref{bir-triv-stab2} of Proposition~\ref{prop: bir+all trivial stabilizers properties}; the only point where the torsion hypothesis is used, is in the implication~\ref{bir-triv-stab3}$\Rightarrow$\ref{bir-triv-stab2}, while the other parts may be repeat ``verbatim''. \\
\ref{characterizationii}$\Rightarrow$\ref{characterizationi}: By Proposition~\ref{prop: positive torsion equivalences2} it is enough to prove that $\St^{+}_{\gauta}(y^{\omega})$, $y\in \XX^{*}$, are formed by torsion elements. Thus, let $\pi(u)\in \St^{+}_{\gauta}(y^{\omega})$ for some~$u\in \QQ^{*}$, $y\in \XX^{*}$. By Proposition~\ref{prop: comm stab}, $u, y^{n}$ is a commuting pair for some integer $n\ge 1$, and by the same proposition $\sigma(y^{n})\in \St_{\gdauta}(u^{\omega})$. Therefore, since $\mathcal{P}(\dual\auta)=\varnothing$, we find $\St^{+}_{\gdauta}(u^{\omega})\neq \{\id\}$, whence $\pi(u^{\ell})=\id$, by Proposition \ref{prop: bir+all trivial stabilizers properties}, \emph{i.e.} $\pi(u)$ is torsion.\\
\ref{characterizationii}$\Rightarrow$\ref{characterizationiii}: Direct by Proposition~\ref{prop: geo description}.\\
\ref{characterizationiii}$\Rightarrow$\ref{characterizationii}: Orbital graphs $\Gamma(\gdauta,\XX, \QQ^{\omega},v^{\omega})$, $v\in \QQ^{*}$ cannot be all finite. Indeed, if $v \in \QQ^*$,  the connected components of its powers have sizes bounded by the size of $\Gamma(\gdauta,\XX, \QQ^{\omega},v^{\omega})$.   Thus, by~\cite[Proposition~7]{KiPiSa15a} $\pi(v)$ is of finite order. And therefore, $\gauta$ is torsion, which contradicts the fact that $\gauta$ has an aperiodic element. Hence there is an infinite (acyclic by hypothesis) orbital graph $\Gamma(\gdauta,\XX, \QQ^{\omega},h^{\omega})$, for some~$h\in \QQ^{*}$. Now, if $\mathcal{P}(\dual\auta) \neq\varnothing$, there is a non-trivial cycle $h^{\omega}\mapright{} h^{\omega}$ in the previous acyclic multidigraph, which is impossible.\\
Finally, let us show (\ref{eq: relations on non-trivial}): in fact, when the orbital graph is finite, then both sides of the equivalence are true, otherwise neither of the sides holds.\\
Let $u\in \QQ^{*}$. If $\Gamma(\gdauta,\XX, \QQ^{\omega},u^{\omega})$ is finite, then clearly $\St_{\gdauta}^{+}(u^{\omega})\neq \{\id\}$ holds and, by the previous argument, $\pi(u)$ is of finite order. On the other hand, if $\Gamma(\gdauta,\XX, \QQ^{\omega},u^{\omega})$ is infinite, then $\pi(u)$ is aperiodic according to~\cite[Proposition~7]{KiPiSa15a}. Further, by hypothesis, the infinite orbital graph $\Gamma(\gdauta,\XX, \QQ^{\omega},u^{\omega})$ is acyclic, whence $\St_{\gdauta}^{+}(u^{\omega})=\{\id\}$, and this concludes the proof.
\end{proof}

In this section we have constructed tools to find singular points, especially in the case of bireversible automata. It is still unknown whether there exists examples of infinite groups generated by RI automata without singular point. In particular we connected singular points with helix graphs and we proved that search can be narrow down to fully positive groups (Theorem~\ref{theo: alternatives}). We also discuss the connection between the existence of non-trivial stabilizer of a group generated by a RI-automaton and the stabilizers of its dual (Theorem~\ref{thm: alternative with dual}). This analysis shows that the possible existence of a RI automaton generating a group with all trivial stabilizers can be restricted to the class of bireversible. However evidences suggest  that if $G$, generated by a RI automaton is infinite then it admits at least a non-trivial boundary. stabilizers.

%-------------------------------------------------------------------------------------------------------------------------------------------------
%-------------------------------------------------------------------------------------------------------------------------------------------------
%-------------------------------------------------------------------------------------------------------------------------------------------------
\section{Dynamics and Wang tilings}\label{sec: wang tilings} 
There is an interesting connection between commuting pairs and Wang tilings observed by~I.~Bondarenko~\cite{Bonda-private}. Given a Mealy automaton, one may associate  a set of Wang tiles reflecting the action of the automaton on the stateset and the alphabet. The existence of periodic tilings corresponds to the notion of commuting words that generate elements of the stabilizers of infinite periodic words. Note that this problem -- called the \emph{domino problem} -- is undecidable in general~\cite{berger,KaPa}. In this section we prove that the domino problem is decidable for some family  of tilesets linked to Mealy automata. On the other hand, we show that the problem of determining whether or not an automaton has commuting words on a restricted stateset is undecidable.

\subsection{Wang tiles \emph{vs} cross-diagrams}
\begin{figure}[h]
\begin{center}
\begin{tikzpicture}[scale=01]
\pavet{0}{0}{$t_n$}{$t_s$}{$t_e$}{$t_w$}
\end{tikzpicture}
\caption{A Wang tile.}
\label{fig: wang tile}
\end{center}
\end{figure}
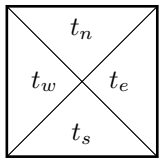

\begin{figure}[h]
\begin{center}
\begin{tikzpicture}[scale=0.8]
\pavet{0}{0}{$v'$}{$v$}{$u'$}{$u$}
\node[] (in) at (3,1)   {$\in \mathcal{T} $};
\node[] (eq) at (4.5,1)   {$\Leftrightarrow $};
\node[] (in) at (7,1)   {{$\begin{array}[b]{ccc}
		& v	&    \\
u		& \lacroix    	& u' \\
		& v'	 			&            		
\end{array}$} };
\node[] (in) at (9,1)   {$\in \auta $};
\node[] (eq1) at (10.5,1)   {$\Leftrightarrow $};
\node[] (aut) at (13,1)   { $u\mapright{v|v'}u' \in \auta$};
\end{tikzpicture} 

\end{center}
\caption{Tileset~$\mathcal{T}(\auta)$ associated with an automaton $\auta$.}\label{fig:associat tiles}
\end{figure}
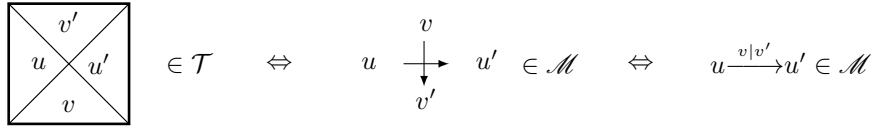

We recall that a Wang tile is a unit square tile with a color on each edge. Formally, it is a quadruple $t=(t_{w}, t_{s}, t_{e}, t_{n})\in C^{4}$ where $C$ is a finite set of colors (see Fig.~\ref{fig: wang tile} for a typical depiction of a Wang tile).
A tileset is a finite set~$\mathcal{T}$ of Wang tiles, and for each $t\in \mathcal{T}$ and~$d\in\{n,s,e,w\}$, we put~$t_{d}$ for the color of the edge in the $d$-side. Given a tileset~$\mathcal{T}$, a tiling of the discrete plane is a map $f: \mathbb{Z}^{2}\rightarrow \mathcal{T}$ that associates to each point in the discrete plane a tile from~$\mathcal{T}$ such that adjacent tiles share the same color on their common edge, \emph{i.e.}, for any $(x,y)\in\mathbb{Z}^{2}$, $f(x,y)_{e}=f(x+1,y)_{w}$ and $f(x,y)_{n}=f(x,y+1)_{s}$. For a rectangle $[x,x']\times [y,y']\subseteq \mathbb{Z}^{2}$, and $d\in \{n,s,e,w\}$ we denote by $f\left([x,x']\times [y,y']\right)_{d}$ the word in $C^{*}$ labelling the $d$-side of the square. For instance for $d=w$ we have:
$$
f\left([x,x']\times [y,y']\right)_{w}=f(x,y)_{w}f(x,y+1)_{w}\cdots f(x,y')_{w}.
$$
This notion extends naturally for rectangles of the form $[x,\infty]\times [y,\infty]\subseteq \mathbb{Z}^{2}$.\\
A tiling $f$ is \emph{periodic} if there exists a periodicity vector
$v\in \mathbb{Z}^{2}$ such that $f(t+v)=f(t)$ for~$t\in\mathbb{Z}^{2}$. A
tiling is \emph{bi-periodic} if there are two linearly independent
vectors~$v_{1}, v_{2}$ for which~$f$ is periodic. It is a well known
fact belonging to the folklore that a tileset admits a periodic tiling
if and only if it admits a bi-periodic tiling if and only if it admits
vertical and horizontal periods. In case  $f:\mathbb{Z}^{2}\rightarrow
\mathcal{T}$ is a periodic tiling with a vertical period $p_{y}$ and a horizontal period~$p_{x}$, a fundamental domain of~$f$ is given by the square $[0,p_{x}-1]\times [0,p_{y}-1]$. Following~\cite{LeGlo14} we say that the tileset~$\mathcal{T}$ is $cd$-deterministic
with~$(c,d)\in\{(e,n),(e,s),(w,n),(w,s)\}$ if each tile~$t\in\mathcal{T}$ is uniquely determined by its pair~$(t_{c},t_{d})$ of colors. Whenever~$\mathcal{T}$ is $cd$-deterministic for each~$(c,d)\in\{(e,n),(e,s),(w,n),(w,s)\}$, we say that $\mathcal{T}$ is \emph{4-way deterministic}.
\\
There is a natural way to associate to a Mealy automaton~$\auta=(\QQ,\XX,\xcdot,\xcirc)$ a tileset~$\mathcal{T}(\auta)$: for each transition~$q\mapright{a|b}p$ we associate the Wang tile~$(q, a, p, b)$ with colors on~$C=\QQ\sqcup \XX$, see Figure~\ref{fig:associat tiles}. This point of view is just a reformulation of the cross-diagram defined previously (and is completely different from the one in~\cite{Gillibert14}).
The following lemma links properties of the automaton~$\auta$ with properties of the associated Wang tileset~$\mathcal{T}(\auta)$.

\begin{lemma}\cite{Bonda-private}
The tileset associated to a  Mealy automaton is necessarily  $ws$-deterministic.
Furthermore, for a Mealy automaton $\auta$ we have the following.

\begin{itemize}
\item$\mathcal{T}(\auta)$ is  $es$-deterministic if and only if $\auta$ is reversible;
\item $\mathcal{T}(\auta)$ is  $wn$-deterministic if and only if $\auta$ is invertible;
\item $\mathcal{T}(\auta)$ is 4-way deterministic if and only if $\auta$ is bireversible.

\end{itemize}
\end{lemma}
\bigskip

\begin{figure}[h] 
\begin{tikzpicture}[scale=0.5]
\pavet{0}{0}{\tiny$v_1' $}{}{}{\tiny$u_k$}
\pavet{0}{-2}{}{}{}{}
\node[]  (vdot1) at (1,-3) {$\vdots$} ;
\pavet{0}{-6}{}{\tiny$v_1$}{}{\tiny$u_1$}

\pavet{2}{0}{\tiny$v_2'$}{}{}{}
\node[]  (vdot22) at (3,-1) {$\vdots$} ;
\node[]  (vdot21) at (3,-3) {$\vdots$} ;
\pavet{2}{-6}{}{\tiny$v_2$}{}{}

\node[]  (vdot33) at (5,1) {$\hdots$} ;
\node[]  (vdot30) at (5,-5) {$\hdots$} ;

\pavet{6}{0}{\tiny$v_\ell'$}{}{\tiny$u_k'$}{}
\node[]  (vdot22) at (7,-1) {$\vdots$} ;
\node[]  (vdot21) at (7,-3) {$\vdots$} ;
\pavet{6}{-6}{}{\tiny$v_\ell$}{\tiny$u_1'$}{}

\node[] (eq) at (9,-2)   {$\Leftrightarrow $};

\cross{11.5}{1}{\tiny v_1}{\tiny u_{1}}{}{}
\node[]  (vdotc1) at (12,-2) {$\vdots$} ;
\cross{11.5}{-5}{}{\tiny u_{k}}{\tiny v_1'}{}
\node[]  (vdotc2) at (14,1) {$\hdots$} ;
\node[]  (vdotc2) at (14,-5) {$\hdots$} ;

\cross{17}{1}{\tiny v_\ell}{}{}{\tiny u_{1}'}
\node[]  (vdotc3) at (16.5,-2) {$\vdots$} ;
\cross{17}{-5}{}{}{\tiny v_\ell'}{\tiny u_{k}'}

\node[] (eq) at (9,-10)   {$\Leftrightarrow    $};
\node[anchor=west] (eq) at (11,-9)   {\footnotesize  $u_1\cdots u_k \xcdot v_1  \cdots v_\ell = v_1'\cdots v_\ell'$};
\node[anchor=west] (eq) at (11,-11)   {\footnotesize $u_1\cdots u_k \xcirc v_1  \cdots v_\ell = u_1'\cdots u_k' $};
\end{tikzpicture}
\caption{Correspondence between partial tilings, cross-diagrams and actions of the automaton.}
\label{fig:wordtil}
\end{figure}
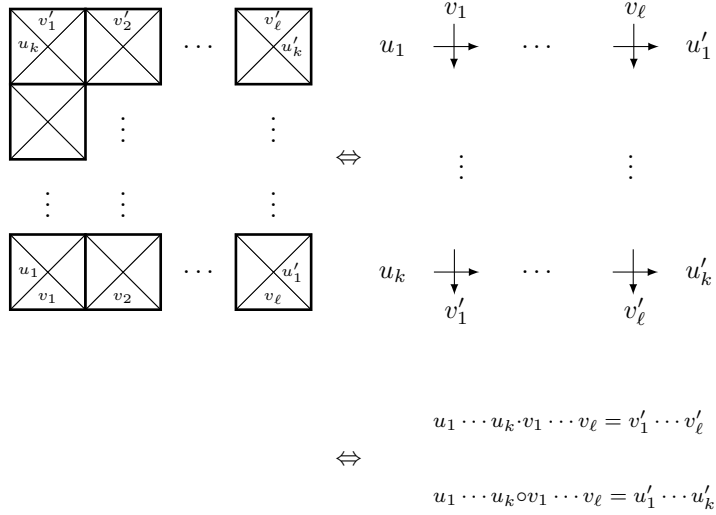

There is a clear correspondence between tilings and cross diagrams, as illustrated in Figure~\ref{fig:wordtil}. This easily extends to infinite words.

The following proposition links periodicity of a tiling and existence of a commuting pair for a Mealy automaton.
\begin{proposition}\label{prop: periodic tiling}
Let $\auta=(\QQ,\XX,\xcdot,\xcirc)$ be an automaton. The following are equivalent.
\begin{enumerate}[label=(\roman{enumi})]
\item $\mathcal{T}(\auta)$ admits a periodic tiling;
\item $\auta$ admits a commuting pair: $v\in \XX^{+}$, $u\in \QQ^{+}$ with $u\xcirc v=v$, $u\xcdot v=u$.
\end{enumerate}
\end{proposition}
\begin{proof}
If $\mathcal{T}(\auta)$ admits a periodic tiling, then there is a bi-periodic tiling $f: \mathbb{Z}^{2}\rightarrow \mathcal{T}(\auta)$
with horizontal and vertical periods $p_{x}, p_{y}$, respectively.
Therefore, the fundamental domain $[0, p_{x}-1]\times [0,p_{y}-1]$ pinpoints a rectangle in the tiling~$f$ with colors given by the string
\[f\left([0, p_{x}-1]\times [0,p_{y}-1]\right)_{n}=f\left([0, p_{x}-1]\times [0,p_{y}-1]\right)_{s}=v\in \XX^{*}\]
in both the south and north edge and
\[f\left([0, p_{x}-1]\times [0,p_{y}-1]\right)_{w}=f\left([0, p_{x}-1]\times [0,p_{y}-1]\right)_{e}=u\in \QQ^{*}\]
in both the west and east edge. Hence, by Figure~\ref{fig:wordtil}: $u\xcirc v=v$, $u\xcdot v=u$.
On the other hand, if there is a commuting pair $u\in \QQ^{*},v\in \XX^{*}$ such that $u\xcirc v=v$,
$u\xcdot v=u$ then by Figure~\ref{fig:wordtil}
there is a partial tiling $f:[0,|v|-1]\times [0,|u|-1]\rightarrow \mathcal{T}(\auta)$ with
\[\begin{array}{c}
f\left([0, |u|-1]\times [0,|v|-1]\right)_{n}=f\left([0, |u|-1]\times [0,|v|-1]\right)_{s}=u,\\
f\left([0, |u|-1]\times [0,|v|-1]\right)_{w}=f\left([0, |u|-1]\times [0,|v|-1]\right)_{e}=v.
\end{array}\]
This tiling extends naturally to a bi-periodic tiling $f:\mathbb{Z}^{2}\rightarrow \mathcal{T}(\auta)$ of the whole discrete plane.
\end{proof}

As a consequence we obtain a decidability result for the domino problem with  particular sets of tilesets.\medskip

Starting from the tileset $\mathcal{T}$ we construct a letter-to-letter transducer $\auta(\mathcal{T})= (\QQ,\XX\xcdot,\xcirc)$ in a natural way as follows: Let $\QQ=\cup_{t\in\mathcal{T}}\{t_{w},t_{e}\}$,  $\XX=\cup_{t\in\mathcal{T}}\{t_{s},t_{n}\}$, and transitions given by $t_w\mapright{t_s|t_n}t_e \in \auta(\mathcal{T}) $ whenever $(t_w,t_s,t_e,t_n)\in\mathcal{T}$. Note that in general $\auta(\mathcal{T})$ is not a Mealy automaton -- we will deal with these cases in the next section. However when it is indeed a Mealy automaton the following holds:

\begin{corollary}\label{cor:period-decid}
Let $\mathcal{T}$ be a Wang tileset. If  $\auta(\mathcal{T})=(\QQ,\XX,\xcdot,\xcirc)$ is a Mealy automaton then $\mathcal{T}$ tiles the plane periodically, in an effective way.
\end{corollary}
\begin{proof}
It is decidable if $\auta(\mathcal{T})$ is a Mealy automaton in time $O( \card{\QQ}\card{\XX} )$. Now construct the helix graph  $\mathcal{H}_{1,1}$ (once again  in time $O( \card{\QQ}\card{\XX} )$).  Since the stateset and the alphabet are finite, and each pair $(u,v) \in \mathcal{H}_{1,1} $ has exactly one successor, there exists a cycle 
\[
(u_{0},v_{0})\mapright{}(u_{1},v_{1})\mapright{}\cdots \mapright{}(u_{m},v_{m})\mapright{}(u_{0},v_{0})
\]
in $\mathcal{H}_{1,1}$. Hence by Lemma~\ref{lem: helix commuting} there is a commuting pair $(u,v) \in \QQ^{+} \times \XX^{+}$, so by Proposition~\ref{prop: periodic tiling}, $\mathcal{T}$ tiles periodically the plane.
\end{proof}
\subsection{Commuting pairs on a restricted stateset}
With Corollary~\ref{cor:period-decid}, we have seen that  each tileset~$\mathcal{T}(\auta)$ associated with a Mealy automaton~$\auta=(\QQ,\XX,\xcdot,\xcirc)$ admits a periodic tiling. However $(u,v)$ could commute because the word $u$ acts like the identity, \emph{i.e.} $\pi(u)=\id$; so this commuting pair will not help us to find singular points, and we may want to avoid it.\\
In this spirit we now consider \emph{commuting pair on a restricted stateset}, that is commuting pair in $\QQ'^*\times \XX^*$, with $\QQ'\varsubsetneq \QQ$. \medskip

We now consider an automaton $\auta= (\QQ,\XX,\xcdot,\xcirc)$ and a subset $\QQ'\varsubsetneq \QQ$.  
We consider the \emph{restricted tileset}~$\mathcal{T}(\auta,\QQ')$ formed by the tiles of~$\mathcal{T}(\auta)$ whose colors are in the set~$\QQ'\cup \XX$. Note that $\mathcal{T}(\auta,\QQ')$ is the tileset associated with the (partial) automaton obtained from~$\auta$ by erasing all the transitions  to or from states in~$\QQ\smallsetminus \QQ'$.
In particular when the automaton has a sink-state\footnote{One can recognize trivial states in linear time (note that we can also find trivial components through \emph{minimization}, with cost~$O(\ |{\XX}||\QQ|\log|\QQ|\ )$, see~\cite{BarSil}). In the following we assume that the automaton has at most one sink-state, which we denote $e$.} $e$,  we define \emph{non-elementary commuting pairs}, that are commuting pairs restricted to the stateset~$\QQ\smallsetminus\{e\}$, and  the \emph{non-elementary tileset}  $\mathcal{T}(\auta,\QQ\smallsetminus\{e\})$ will be denoted by $\overline{\mathcal{T}}(\auta)$. We have the following proposition.

\begin{proposition}\label{prop: reduced tileset}
Let  $\auta $ be a Mealy automaton with a sink-state. Then $\auta$ admits a non-elementary commuting pair if and only if $\overline{\mathcal{T}}(\auta)$ has a periodic tiling.
\end{proposition}
\begin{proof}
The proof is similar to the one presented in Proposition~\ref{prop: periodic tiling}. Indeed if $e$ were to appear in the cross diagramm, then since it a sink it will stay present  every step after, preventing the tiling to by periodic: $\forall u \in \QQ^*e\QQ^*, \forall v \in \XX^*, u \xcdot v \in \QQ^*e\QQ^*$. 
  
\end{proof}

This raises the problem of finding restricted commuting pairs, since we can no longer apply Corollary~\ref{cor:period-decid}. We obtain the following decision problem:\medskip

\newcommand\RCP{\small\texttt{RESTRICTED\! COMMUTING\! PAIRS}}

{\RCP:}
\begin{itemize}
\item \textbf{Input:} $\auta=(\QQ,\XX,\xcdot,\xcirc)$, $\QQ'\varsubsetneq \QQ$.
\item \textbf{Output:} Does $\auta$ have commuting pairs restricted to the stateset $\QQ'$?
\end{itemize}

\newcommand\NECP{\small\texttt{NON-ELEMENTARY\! COMMUTING\! PAIRS}}

{It turns out that this problem is undecidable. To prove this we are going to reduce a  known undecidable problem to it, namely the existence of a periodic tiling for a 4-way deterministic tileset (see \cite{LeGlo14}).\medskip

For the reduction from the periodic 4-way deterministic problem to our problem we show that, given a 4-way deterministic tileset, a periodic tiling can be determined from a non-trivial commuting pair of some Mealy automaton.\\

\begin{figure}[h] \label{fig:tilltoaut}
\begin{tikzpicture}%
\begin{scope}[scale=0.48]
\pavet{0}{0.5}{\tiny a}{\tiny b}{\tiny 1}{\tiny 1}
\pavet{3}{0.5}{\tiny a}{\tiny b}{\tiny 3}{\tiny 2}
\pavet{0}{-2.5}{\tiny b}{\tiny a}{\tiny 1}{\tiny 3}
\pavet{3}{-2.5}{\tiny b}{\tiny a}{\tiny 2}{\tiny 2}
\end{scope}

\node[] (eq) at (3.3,0)   {$\Longrightarrow $};

\begin{scope}[xshift=120]
	\node[state] (a) at (0,0)   {\tiny $a$};
	\node[state] (b) at (2,0)  {\tiny $b$};
	\path
		(a)	edge[bend left,->]	node[above]{\tiny \begin{tabular}{c}\(1|2\)\\\(2|3\)\\\end{tabular}}			(b)
		(b)	edge[bend left,->]	node[below]{{\tiny \begin{tabular}{c}\(2|2\)\\\(3|1\)\\\end{tabular}}}	(a);
\end{scope}	

\node[] (eq) at (8,0)   {$\Longrightarrow    $};

\begin{scope}[xshift=260]
	\node[state] (a') at (0,0.5)   {\tiny $a$};
	\node[state] (b') at (2,0.5)  {\tiny $b$};
	\node[state] (e) at (1,-1)  {\tiny $e$};
	\path
		(a')	edge[bend left,->]	node[above]{\tiny \begin{tabular}{c}\(1|2\)\\\(2|3\)\\\end{tabular}}			(b')
		(b')	edge[bend left,->]	node[below]{{\tiny \begin{tabular}{c}\(2|2\)\\\(3|1\)\\\end{tabular}}}	(a')
		(a')	edge[->]			node[left]{\tiny \(3|1\)}			(e)
		(b')	edge[->]			node[right]{\tiny \(1|3\)}			(e)
		(e)	edge[loop left] 		node[left=-4]{\tiny \begin{tabular}{c}\(1|1\)\\\(2|2\)\\ \(3|3\)\end{tabular}}	(e);
\end{scope}
\end{tikzpicture}
\caption{Construction of an invertible  Mealy automaton from a 4-way deterministic tileset.}
\label{fig-4way}
\end{figure}
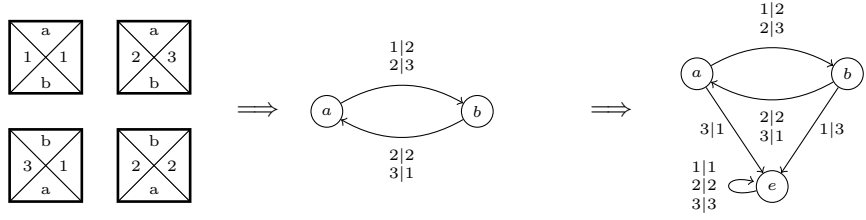

Let $\mathcal{T}$  be a 4-way deterministic tileset. We build $\auta(\mathcal{T})$ as in Corollary~\ref{cor:period-decid}.\\
This transducer is not in general complete, however, we may add an extra sink-state $e$ and extra transitions to this sink-state in order to make the automaton complete and invertible. We call $\auta_{\mathcal{S}}(\mathcal{T}) $ the invertible automaton associated with~$\mathcal{T}$. Since the extra transitions go to the sink-state it is clear that $\mathcal{T}=\overline{\mathcal{T}}(\auta_{\mathcal{S}}(\mathcal{T}))$.\\
We are now in position to prove the following result.

\begin{theorem}\label{theo: undecidability non-elementary comm}
The problem {\rm\RCP} is  undecidable.
\end{theorem}
\begin{proof}
We prove that a particular instance of $\RCP$ is undecidable: \\
{\NECP:}
\begin{itemize}
\item \textbf{Input:} $\auta$ a Mealy automaton with a sink-state.
\item \textbf{Output:} Does $\auta$ have non-elementary commuting pairs?
\end{itemize}
By~\cite{LeGlo14} checking whether a 4-way deterministic tileset admits a periodic tiling is undecidable. By Proposition~\ref{prop: reduced tileset} checking whether $\mathcal{T}$ admits a periodic tiling is equivalent to check whether  $\auta_{\mathcal{S}}(\mathcal{T})$ has a non-elementary commuting pair, hence {\rm\NECP} is undecidable, and so is \RCP.
\end{proof}
Note that in this context we deal with inverses of~$\QQ$, and it is not difficult to check that reduced words $u\in \wt{\QQ}^{*}$ may generate non-elementary commuting pair. Hence, in this context a  non-elementary pair of commuting words is a pair of words $u\in \wt{\QQ\smallsetminus\{e\}}^{*}, v\in \XX^{*}$ such that $u$ is reduced and $u,v$ commute.

\begin{proposition}\label{prop: not continuous with inverses}
With the above notation. Let $\auta\in\SinkA$. The following are equivalent.
\begin{enumerate}[label=(\roman{enumi})]
\item\label{not-cont-inv1}$v^{\omega}$ is a singular point, $v\in \XX^{*}$;
\item\label{not-cont-inv2} there is a non-elementary commuting pair $u\in \wt{\QQ\smallsetminus\{e\}}^{*}$, $v\in \XX^{*}$ with \mbox{$\pi(u)\neq \id$}.
\end{enumerate}
\end{proposition}
\begin{proof}
\ref{not-cont-inv1}$\Rightarrow$\ref{not-cont-inv2} If $v^{\omega} \in \krit$ then the map $\SSt$ is not continuous at~$v^{\omega}$ if
  and only if $\St_{\gauta}(v^{\omega})\neq \St_{\gauta}^{0}(v^{\omega})$. This
  implies that there exists $u \in \wt{\QQ\smallsetminus\{e\}}^{*}$ satisfying~$\pi(u)\in
  \St_{\gauta}(v^{\omega})\smallsetminus\St_{\gauta}^{0}(v^{\omega})$. Note that $\pi(u\xcdot v^k)\neq \id$ for any $k\geq 0$ (from Lemma~\ref{lemmarestric}). Since $|u|$ is finite, there exist $n>m>0$ such that $u\xcdot v^m=u\xcdot v^n$. Then $u_m:=u\xcdot v^m$ and $v^{n-m}$ form  a non-elementary commuting pair with~$u_{m}\in \wt{\QQ\smallsetminus\{e\}}^{*}$. \\
\ref{not-cont-inv2}$\Rightarrow$\ref{not-cont-inv1} Let $u\in \wt{\QQ\smallsetminus\{e\}}^{*}$, $v\in
  \XX^{*}$ be a non-elementary commuting pair. Since $\pi(u)\neq \id$, there
  is a word~$w\in \XX^{*}$ such that $u\xcirc w\neq w$. Consider the
  sequence $w_n:=v^nwv^{\omega}$ and proceed as in the proof of
  Lemma~\ref{lemmaperiodiche} to show that $\pi(u)\in
  \St_{\gauta}(v^{\omega})\smallsetminus \St_{\gauta}^{0}(v^{\omega})$.
\end{proof}

\begin{rem}
Note that Proposition \ref{prop: discontinuous point birev} implies that Proposition~\ref{prop: not continuous with inverses} holds also for bireversible automata. 
\end{rem}

As we have already noted, the tileset~$\mathcal{T}(\auta)$
associated with a Mealy automaton~$\auta$ always admits a periodic
tiling, however, in case we consider automata from~$\SinkA$, it
is interesting, and useful for the sequel, to understand when the non-elementary tileset of~$\overline{\mathcal{T}}(\auta)$  tiles the discrete plane. We need first the following lemma regarding inverse $X$-digraphs, that are digraph $\Gamma$ with edges labelled by element of $X$,  such that if $p\mapright{a}q$ is an edge of $\Gamma$, so is $q\mapright{a^{-1}}p$ (involutive), for any vertex $v\in \Gamma $ and any label $x \in \wt{X}$, there is exactly one edge starting from $v$ and labelled with $x$. In inverse digraph the natural quasimetric $d$ is symmetric and so it is a distance.\\
Note that an inverse $X$-digraph has in- and out-degree $\card{\wt{X}}$. 

\begin{lemma}\label{lem: infinite path}
Let $\Gamma$ be an infinite connected inverse $X$-digraph, $|X|<\infty$, and let $v$ be a vertex of \/ $\Gamma$. Then there is a right-infinite word $\theta\in \wt{X}^{\omega}$ such that \[\lim_{i\rightarrow \infty}d(v, v_{i})=\infty \: ,\] where $v_i$ is defined by {$v\longmapright{\theta[:i]} v_{i}$}.
\end{lemma}

\begin{proof}
Since $\Gamma$ is an inverse digraph it is symmetric. Consider $\overline{\Gamma}$ the non-directed graph obtained by gluing edges $p\mapright{a}q$ and $q\mapright{a^{-1}}p$ in $\Gamma$. Note that the distances in $\Gamma$ and in $\overline{\Gamma}$ are the same and we will denote both of them by~$d(.,.)$. Since $\overline{\Gamma}$ has finite degree, is infinite and connected, we can apply K\"onig's lemma that claims that $\overline{\Gamma}$ contains an infinitely long simple path (that is, a path with no repeated vertices). Let $\nu = v_1\cdots v_i \cdots$ denote such a path and $\theta$ the associated infinite word. Since $\overline{\Gamma}$ is connected there exists a path from any vertex $v$ to $v_1$: let $k=d(v,v_1)$. Since the graph~$\overline{\Gamma}$ has bounded degree, the ball~$\ball_k(v)$ of radius $k$ centred in $v$  is finite and since $\nu$ is simple there is an index $i_k$ for which $v_j \notin\ball_k(v)$ for~$j\geq i_k$. The same holds for arbitrary large~$k$, which concludes the proof.
\end{proof}

Consider the map $\lambda$ from the set of edges of~$\Sch(\xi)$ to the integers, that associate to each edge~$f=\eta\mapright{q}\eta'$ the integer
$$
\lambda(f)=
\begin{cases} \min\{n: q\xcdot \eta[:n]=e\} &\mbox{if } \exists m:  q\xcdot \eta[:m]=e\\
\infty & \mbox{otherwise}.
\end{cases}
$$
Let $d(\eta,\eta')$ denote the usual metric on $\Sch(\xi)$. Since
$\Sch(\xi)$ is a regular digraph with finite out-degree $\wt{\QQ}$
for each integer $n$ we consider the set of edges that admit at least one vertex inside the ball of radius $n$ centered at~$\xi$:
$$
\balle_{\xi}(n)=\left\{\eta\mapright{q}\eta': d(\xi, \eta)\le n\vee d(\xi, \eta')\le n\right\}.
$$
Note that $\balle_{\xi}(n)$ is clearly finite and we may define
\[\psi_{\xi}(n)=\max\left\{\lambda(f): f\in \balle_{\xi}(n)\wedge \lambda(f)<\infty\right\}.\]
We remark that $\left\{\lambda(f): f\in \balle_{\xi}(n)\wedge \lambda(f)<\infty\right\}$ is non-empty in view of the existence of the trivial edges of type $\eta\mapright{e}\eta$. Thus, $\psi_{\xi}$ is a monotonically increasing function. We have thus two cases: either $\lim_{n\rightarrow \infty} \psi_{\xi}(n)=\infty$, or $\lim_{n\rightarrow \infty} \psi_{\xi}(n)=\ell$ for some~$\ell\ge 1$. The following proposition provides sufficient conditions on the dynamics of the boundary for the non-elementary tileset to tile the discrete plane.

\begin{proposition}\label{prop: almost equiv conditions to have a tiling}
Let $\auta=(\QQ,\XX,\xcdot,\xcirc)\in\SinkA$. If there is a point~$\xi\in \XX^{\omega}$ whose Schreier graph $\Sch(\xi)$ is infinite and such that $\lim_{n\rightarrow \infty} \psi_{\xi}(n)=\ell$ for some integer~$\ell\ge 1$, then $\overline{\mathcal{T}}(\auta)$ tiles the discrete plane. Conversely, if $\overline{\mathcal{T}}(\auta)$ admits just aperiodic tilings, then there is an infinite Schreier graph $\Sch(\xi)$ for some~$\xi\in \XX^{\omega}$. 
\end{proposition}
\begin{proof}
It belongs to the folklore that a tileset tiles the whole discrete plane
if and only if it tiles the first quadrant of the discrete plane
(see for instance~\cite{Robinson}).
Let us prove that there is a tiling $f:[1,\infty]\times[1,\infty]\rightarrow \overline{\mathcal{T}}(\auta)$ of the first quadrant. First there is a tiling   $h:[1,\infty]\times[1,\infty]\rightarrow \mathcal{T}(\auta)$ of the first quadrant such that, putting
$
\eta=h\left([1,\infty]\times [1,\infty]\right)_{s}$ and $  \sigma=h\left([1,\infty]\times [1,\infty]\right)_{w}
$, we have:
\begin{enumerate}
\item \label{prop-aper-proof1}$\sigma$ is not cofinal with~$e^{\omega}$, and
\item \label{prop-aper-proof2}if $\sigma[j]\neq
e$, for some~$j\ge 1$, then $\sigma[j]\xcdot (\sigma[:j-1]\xcirc\eta[:i])\neq e$ for all $i\ge 1$.
\end{enumerate}
\begin{figure}[h] \label{fig:tillschreier}
\begin{tikzpicture}%
\def\cori{1}
\def\corj{4}
\def\corjj{2}

\coordinate (A) at (0*\corj,0*\cori); 			\coordinate (Aw) at (-1*\corjj,0*\cori); 	\coordinate (Aww) at (-2*\corjj,0*\cori);
\coordinate (B) at (0*\corj,-1*\cori); 			\coordinate (Bw) at (-1*\corjj,-1*\cori); 	\coordinate (Bww) at (-2*\corjj,-1*\cori);
\coordinate (C) at (0*\corj,-2*\cori);
\coordinate (D) at (0*\corj,-3*\cori);			\coordinate (Dw) at (-1*\corjj,-3*\cori);	\coordinate (Dww) at (-2*\corjj,-3*\cori);
\coordinate (E) at (0*\corj,-4*\cori); 			\coordinate (Ew) at (-1*\corjj,-4*\cori); 	\coordinate (Eww) at (-2*\corjj,-4*\cori);
\coordinate (F) at (0*\corj,-4.5*\cori);

\coordinate (A') at (1*\corj+0*\corj,0*\cori); 			\coordinate (Aw') at (1*\corj+1*\corjj,0*\cori); 	\coordinate (Aww') at (1*\corj+2*\corjj,0*\cori);
\coordinate (B') at (1*\corj+0*\corj,-1*\cori); 			\coordinate (Bw') at (1*\corj+2*\corjj,-1*\cori); 	\coordinate (Bww') at (-2*\corjj,1*\cori);
\coordinate (C') at (1*\corj+0*\corj,-2*\cori);
\coordinate (D') at (1*\corj+0*\corj,-3*\cori);			\coordinate (Dw') at (1*\corj+1*\corjj,-3*\cori);	\coordinate (Dww') at (1*\corj+2*\corjj,-3*\cori);										
\coordinate (E') at (1*\corj+0*\corj,-4*\cori); 			\coordinate (Ew') at (1*\corj+1*\corjj,-4*\cori); 	\coordinate (Eww') at (1*\corj+2*\corjj,-4*\cori);
\coordinate (F') at (1*\corj+0*\corj,-4.5*\cori); 	

\coordinate (A'') at (2*\corj+0*\corj,0*\cori); 			\coordinate (Aw'') at (2*\corj+1*\corjj,0*\cori); 	\coordinate (Aww'') at (2*\corj+2*\corjj,0*\cori);
\coordinate (B'') at (2*\corj+0*\corj,-1*\cori); 			\coordinate (Bw'') at (2*\corj+2*\corjj,-1*\cori); 	\coordinate (Bww'') at (-2*\corjj,1*\cori);
\coordinate (C'') at (2*\corj+0*\corj,-2*\cori);
\coordinate (D'') at (2*\corj+0*\corj,-3*\cori);			\coordinate (Dw'') at (2*\corj+1*\corjj,-3*\cori);	\coordinate (Dww'') at (2*\corj+2*\corjj,-3*\cori);									
\coordinate (E'') at (2*\corj+0*\corj,-4*\cori); 			\coordinate (Ew'') at (2*\corj+1*\corjj,-4*\cori); 	\coordinate (Eww'') at (2*\corj+2*\corjj,-4*\cori);

\draw[] (E) --  (E'') node[above,near start] {\(~{\xi}\)};

\draw[very thick,->] (E') --  (A') node[left, midway] {\(~{\sigma}\)};
\draw[very thick,->] (E') --  (E'') node[below, midway] {\(~{\eta}\)};

\draw[] (A) --  (E) node[left,midway] {\(~{\theta}\)};

\draw[decorate,decoration=brace] (F') -- (F) node[below,midway] {\(~\ell\)};

\draw[pattern=north west lines, pattern color=blue] (E') rectangle (7.9*\cori,-0.1*\corjj);

\node[rectangle,draw, thin, fill=white]  at (5.5*\cori,-1.5*\corjj) {\tiny{words $e^{\omega}$} or $e$ free} ;
\draw[densely dotted,thick] (B'') --  (0.33*\corj+0*\corj,-1*\cori) node[left] {\(~\tiny{e}\)};

\end{tikzpicture}
\caption{Construction of a non-elementary tilling when $\lim_{n\rightarrow \infty} \psi_{\xi}(n)=\ell$.}
\label{fig:nonelemtil}
\end{figure}
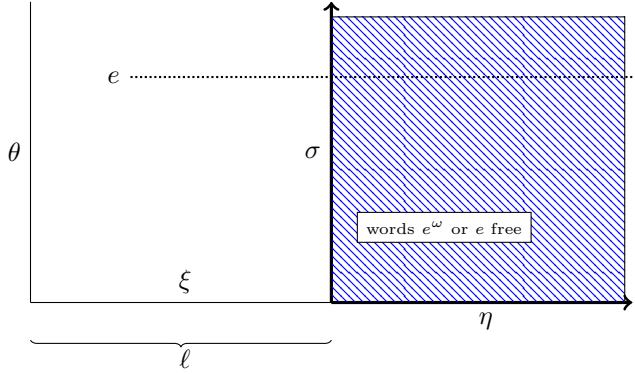

 Intuitively, the last property means that if the west edge of the leftmost tile of a band is not colored by $e$, then also none of the west and east edges of all in this band are colored by $e$. Indeed, since $e$ is a sink, if the west edge of the leftmost tile of a band is colored by $e$, then clearly all the west and east edges of the all the tiles in this band are colored by $e$ as well, see Fig.~\ref{fig:nonelemtil}. Let us first find a tiling fulfilling Condition \ref{prop-aper-proof2}:
If there exists an element $q\in \QQ$ of infinite order, take~$q^\omega$ as $\sigma$ and some word of infinite orbit under the action of $q$ as $\eta$. Otherwise, since $\Sch(\xi)$ is infinite we may find, by Lemma~\ref{lem: infinite path}, a right-infinite path starting at~$\xi$ and labelled by some right-infinite word $\theta\in \wt\QQ^{\omega}$ such that $\lim_{i\rightarrow \infty}d(\xi, \theta[:i]\xcirc \xi)=\infty$. Since every element of $\QQ$ has finite order we may remplace each negative occurence of a letter by the suitable power of this letter, hence without loss of generality we can suppose $\theta \in \QQ^\omega$. By Figure~\ref{fig:wordtil} consider the  tiling $h:[1,\infty]\times[1,\infty]\rightarrow \mathcal{T}(\auta)$ associated with~$\xi, \theta$, \emph{i.e.}, such that
$$
\xi=h\left([1,\infty]\times [1,\infty]\right)_{s}, \; \theta=h\left([1,\infty]\times [1,\infty]\right)_{w}\/.
$$
Since $\lim_{n\rightarrow \infty} \psi_{\xi}(n)=\ell$, then for any vertex $\xi'$ of~$\Sch(\xi)$ and $q\in \QQ$ the following property:
\begin{equation}\label{eq: property tiling}
\mbox{ if }q\xcdot \xi'[:i]\neq e \mbox{ for all }1\le i\le \ell, \mbox{ then } q\xcdot \xi'[:i]\neq e\mbox{ for all }i\ge \ell \phantom{\:.}
\end{equation}
holds. Since $e$ is a sink-state, Equation~\eqref{eq: property tiling} can be written equivalently:  
\begin{equation}
\tag{\ref{eq: property tiling}}
\text{if } q \xcdot \xi'[:\ell] \neq e, \phantom{\mbox{ for all }1\le i\le \ell} \mbox{ then } q\xcdot \xi'[:i]\neq e\mbox{ for all }i\ge \ell
\:.\end{equation}
Let us restrict the tiling $h$ to the quadrant ${{[\ell, \infty]}\times [1,\infty]}$. Take
\begin{align*}
\eta &= h({{[\ell, \infty]}\times [1,\infty]})_s = \xi[\ell :] \:, \\
\sigma&= h({{[\ell, \infty]}\times [1,\infty]})_w = \theta \xcdot  \xi[:\ell-1] \:.
\end{align*}
Property~\ref{prop-aper-proof2} is obtained from Equation~\eqref{eq: property tiling}. \\
Let us prove now Property~\ref{prop-aper-proof1}. Two elements~$\xi', \xi''\in \XX^{\omega}$ are $m$-cofinal if $\xi'[m:]=\xi''[m:]$. Note that for fixed $\xi'\in \XX^{\omega}$ and $m\ge 1$, the number of vertices $\xi''$ that are $m$-cofinal with
\(\xi'\) is finite. Suppose that $\sigma$ is  cofinal with~$e^{\omega}$: $\sigma=\sigma[:k]e^{\omega}$ for some~$k\ge 1$. Hence  the sequence $(\sigma[:i]\xcirc \eta)_i$ is ultimately constant. But $\sigma[:i]\xcirc \eta$ is the suffix of  $ \theta[:i] \xcirc \xi$, whence the set $\{ \theta[:k+j] \xcirc \xi: j\ge 1\}$  formed  $\ell$-cofinal vertices is finite. Which contradicts $\lim_{i\rightarrow \infty}d(\xi,  \theta[:i] \xcirc \xi)=\infty$, and $\sigma$ cannot be cofinite with $e^\omega$. Whence we have obtained a tiling of the first quadrant with $\sigma$ and $\eta$ as west and south borders respectively, and such that, for all $i$, either $\sigma[i] \xcdot \eta  = e^\omega$ or $\sigma[i] \xcdot \eta$ contains no $e$. We now obtain a tiling in $\overline{\mathcal{T}}(\auta)$  by deleting those line $e^\omega$.\medskip

Let us prove the converse in our proposition. Let 
$$
\xi=f\left([1,\infty]\times [1,\infty]\right)_{s}, \;\theta=f\left([1,\infty]\times [1,\infty]\right)_{w}
$$
where $f: \mathbb{Z}\times\mathbb{Z}\rightarrow \overline{\mathcal{T}}(\auta)$ is an aperiodic tiling. Since $f$ is aperiodic, we get $\theta[:i]\xcirc \xi\neq \theta[:j]\xcirc \xi$ for $i\neq j$ (since periodicity in one direction is equivalent to periodicity in both directions).  Therefore, $\{\theta[:i]\xcirc \xi: i\ge 0 \}$ is infinite and $\Sch(\xi)$ as well.

\end{proof}

\newcommand\MaxSync{\small\textbf{{Maximal~non~Synchronizing }}}

We recall that an automaton $\auta=(\QQ,\XX,\xcdot)$ is \emph{synchronizing} whenever there is a word $w\in \XX^{*}$ such that $q\xcdot w=p\xcdot w$ for any $q,p\in \QQ$.
A Mealy automaton~$\auta = (\QQ,\XX,\xcdot,\xcirc)$ is synchronizing whenever $(\QQ,\XX,\xcdot)$ is synchronizing. The set of synchronizing words is denoted by~$\Syn(\auta)$.
Note that an automaton in~$\SinkA$ is always synchronizing. For further details on synchronizing automata and some connections with automaton groups see~\cite{DaRo13}. The next proposition characterizes automata in $\SinkA$ whose set of reduced tiles does not tile the discrete plane in terms of a stability property regarding synchronization. \\
Consider the prefix-closed  language \[\NSyn_{\QQ}(v)=\{u\in \QQ^{<{|v|-1}} \mid \forall j \le |u|, \: u[:j] \xcirc v\notin\Syn(\auta)\}\:.\] and the property~{\MaxSync}: there exists an integer $m$ such that, for all $v \in \XX^{m}$, any  $u$  maximum in $\NSyn_{\QQ}(v)$ for the prefix relation satisfies $uq \xcirc v\in \Syn(\auta)$, for all $q\in \QQ$.

\begin{proposition}\label{prop: no tiling characterization}
Let $\auta=(\QQ,\XX,\xcdot,\xcirc)\in\SinkA$. Then $\overline{\mathcal{T}}(\auta)$ does not tile the discrete plane if and only if property~{\MaxSync} holds for $\auta$.
\end{proposition}
\begin{proof}
Suppose that $\overline{\mathcal{T}}(\auta)$ does not tile the discrete plane. Hence there exists an integer $m$ such that $\overline{\mathcal{T}}(\auta)$ cannot tile a square of size $m$. 
\\
Let $v \in \XX^{m}$ and $u \in \NSyn_{\QQ}(v)$. First suppose $|u|<m-2$: if for some $q \in \QQ$ we have $uq\xcirc v\notin \Syn(\auta)$ then $uq \in \NSyn_{\QQ}(v)$, hence $u$ cannot be maximal. Now suppose $|u|=m-2$:
 if $uq \xcirc v\notin \Syn(\auta)$, for some~$q\in \QQ$, then there is a $q'\in \QQ$ such that $q'\xcdot (uq\xcirc v)[:j]\neq e$ for all $j\le m$. Hence, by Figure~\ref{fig:wordtil} there is a tiling of the square of size $m$ associated with the two words $v,uqq'$, a contradiction.
\\
Conversely, if $\overline{\mathcal{T}}(\auta)$ tiles a square $f:[1,m]\times [1,m]\rightarrow \overline{\mathcal{T}}(\auta)$, then if we put
$$
h=f\left([1,m]\times [1,m]\right)_{s}, \;v=f\left([1,m]\times [1,m]\right)_{w}
$$
we get that $v[:j]\xcirc h\notin \Syn(\auta)$ for all $ j\le m-1$. In particular,  $u=v[:m-2]$  is maximal in $\NSyn_{\QQ}(v)$. However, $(u v[m-1])\xcirc h=v[:m-1]\xcirc h \notin \Syn(\auta)$, hence the condition of the statement is not satisfied.
\end{proof}

\subsection{Aperiodic tilings and singular points}

In Subsection~\ref{subsec: ri-automata} we have seen that $RI$-automata, with empty set of singular points are necessarily bireversible with all trivial stabilizers in the boundary (see Corollary~\ref{cor: ri-automaton continuous point are bireversible} and Proposition~\ref{prop: discontinuous point birev}). Furthermore, Corollary~\ref{cor: empty positive relations implies stab} implies that if a bireversible automaton~$\auta$ has no positive relations, then it must have a non-trivial stabilizer in the boundary. In particular, for the class of bireversible automata it is not possible to have simultaneously no positive relations and no singular points. \medskip

This fact no longer holds in the class~$\SinkA$: there exist automata in $\SinkA$ without positive relations and with no singular points. However, we need some precaution in defining the set of
``positive relations'' for an automaton~$\auta\in\SinkA$.
Indeed, for such an automaton, the set 
${\mathcal{P}}(\auta)$ is always non-empty since the sink-state~$e$ acts
like the identity. 
\\
In order to characterize tilings without non-trivial commuting pairs, we need to define a special set of relations that are in some sense non elementary. Such relations intuitively correspond to words with the property that some of their restrictions do not become the trivial word. For $w\in\QQ^{*} $, $|w|_{e}$ denotes the number of occurrences of the letter $e$ in $w$. The set $\mathcal{E}(\auta)$ of the \emph{non-elementary relations} is defined by:
\[\mathcal{E}(\auta)=\left\{u\in (\QQ\smallsetminus\{e\})^{*}: \pi(u)=\id, \;\exists v\in\XX^{*}, |u\xcdot v^{n}|_{e}<|u|, \,\forall n\ge 1\right\}.\]
By a compactness argument note that the complement of the non-elementary relations, the set of \emph{elementary relations}, may be described as the set of words~$u\in (\QQ\smallsetminus\{e\})^{*}$ with $\pi(u)=\id$, such that there exists some~$n\ge 1$ for which $u\xcdot v=e^{|v|}$ for every~$v\in \XX^{\ge n}$.
Geometrically, automata in the class~$\SinkA$ with no singular points and whose eventual relations are elementary relations, possess helix graphs with a particular shape, as proved in Proposition~\ref{prop: equivalence positive continuity}. We
say that the helix graph~$\mathcal{H}_{k,n}$ of an automaton~$\auta\in\SinkA$ is \emph{singular} whenever each
connected component of~$\mathcal{H}_{k,n}$ has a unique cycle which is necessarily of the
form~$(e^{k}, v)\mapright{}(e^{k},v)$ for some~$v\in \XX^{n}$.

\begin{proposition}\label{prop: equivalence positive continuity}
Let $\auta\in\SinkA$. The following are equivalent:
\begin{enumerate}[label=(\roman{enumi})]
\item \label{equ-pos-cont1}for any $\xi\in \XX^{\omega}$, $g\in\St^{+}_{\gauta\!\!}(\xi)$, there
  is $n\ge 1$ with $g\xcdot \xi[:\!n]=\id$ and~$\mathcal{E}(\auta)=\varnothing$;
\item \label{equ-pos-cont2}for any $v\in \XX^{*}$, $g\in\St^{+}_{\gauta\!\!}(v^{\omega})$, there
  is $n\ge 1$ with~$g\xcdot v^n=\id$ and~$\mathcal{E}(\auta)=\varnothing$; 
\item \label{equ-pos-cont3}for any~$k,n\ge 1$, the helix graph $\mathcal{H}_{k,n}$ is singular;
\item \label{equ-pos-cont4}there is no non-elementary pair of commuting words.
\end{enumerate}
\end{proposition}
\begin{proof}
\ref{equ-pos-cont1}$\Rightarrow$\ref{equ-pos-cont2}. Trivial. \\
\ref{equ-pos-cont2}$\Rightarrow$\ref{equ-pos-cont3}. Suppose that $\mathcal{H}_{k,n}$ is not
  singular for some~$k,n\ge 1$. Therefore, by Lemma~\ref{lem: helix
    commuting} there is a commuting pair $u\in
  \QQ^{*}\smallsetminus\{e\}^{*}, v\in \XX^{*}$. As $e$ acts like the identity, by erasing the (potential) occurrences of $e$ in $u$, we obtain a word  $u'\in (\QQ\smallsetminus\{e\})^{*}$ such that $u', v$ commutes. Hence Lemma~\ref{lem: helix commuting} implies $\pi(u')\in \St^{+}_{G}(v^{\omega})$. If $\pi(u')=\id$, then, since $u',v$ commutes, $u'\in \mathcal{E}(\auta)\neq \varnothing$. Otherwise, we have $\pi(u')\xcdot v^{n}=\pi(u')\neq \id$ for all $n\ge 1$. \\
\ref{equ-pos-cont3}$\Rightarrow$\ref{equ-pos-cont4}. If there is a non-elementary commuting pair $u\in (\QQ\smallsetminus\{e\})^{*}, v\in \XX^{*}$, then by Lemma~\ref{lem: helix commuting} the helix graph~$\mathcal{H}_{|u|,|v|}$ contains the loop~$(u,v)\mapright{} (u,v)$, \emph{i.e.}, $\mathcal{H}_{|u|,|v|}$ is not singular. \\
\ref{equ-pos-cont4}$\Rightarrow$\ref{equ-pos-cont1}. If $\mathcal{E}(\auta)\neq\varnothing$, then for any $u\in\mathcal{E}(\auta)$, by definition, there exists $v\in \XX^{*}$ satisfying
\[|u\xcdot v^{n}|_{e}<|u|, \,\forall n\ge 1.\]
By a compactness argument, there exists some~$m\ge 1$ with~$|u\xcdot v^{m}|_{e}=|u\xcdot v^{m+k}|_{e}$ for every~$k\ge 0$. Moreover, there exist indices~$i>j\ge m$ satisfying~$u\xcdot v^{j}=u\xcdot v^{i}$. As $e$ acts like the identity, by erasing the (potential) occurrences of $e$ in $u$, we obtain a word  $u'\in (\QQ\smallsetminus\{e\})^{*}$ such that $u', v^{i-j}$ is a non-elementary pair of commuting words. Thus, we may assume $\mathcal{E}(\auta)=\varnothing$. Now suppose~$g\in\St^{+}_{\gauta}(\xi)$ for some~$\xi\in \XX^{\omega}$ with~$g\xcdot \xi[:n]\neq \id$ for all~$n\ge 1$. Hence, there are some word~$u\in (\QQ\smallsetminus\{e\})^{*}$ such that $\pi(u)=g$ and some indices~$i> j\ge 1$ such that if we put
$u'=u\xcdot \xi[:j]\in (\QQ\smallsetminus\{e\})^{*}$, then $u'$, $\xi[j+1:i]$ is a non-elementary commuting pair.
\end{proof}

Note that the previous proposition provides necessary conditions on~$\auta\in\SinkA$ for tilings of~$\overline{\mathcal{T}}(\auta)$ of the discrete plane to be aperiodic. As a result of Theorem~\ref{theo: undecidability non-elementary comm} we immediately obtain the following undecidability result of checking the previous ``continuity'' condition.

\begin{theorem}
Given an automaton~$\auta\in\SinkA$, it is undecidable whether for any~$\xi\in \XX^{\omega}$, $g\in\St^{+}_{G}(\xi)$, there exists~$n\ge 1$ with~$g\xcdot \xi[:n]=\id$ and~$\mathcal{E}(\auta)=\varnothing$.
\end{theorem}

Moreover, by Proposition~\ref{prop: equivalence positive continuity},
taking the aperiodic $4$-way deterministic tileset~$\mathcal{T}$
described in~\cite{KaPa} and the associated 
automaton~$\auta\in\SinkA$ with~$\mathcal{T}=\overline{\mathcal{T}}(\auta)$
we get that $\auta$ actually satisfies the ``continuity'' conditions
described in Proposition~\ref{prop: equivalence positive
  continuity}. The same paper raised the problem to determining
the existence of an aperiodic reflection-closed tileset.
A tileset~$\mathcal{T}$ with colored oriented edges is \emph{closed under reflection}
if for each tile in~$\mathcal{T}$ the reflection of this tile along a
horizontal or vertical line also belongs to~$\mathcal{T}$. Kari and
Papasoglu considered the following tiling rule: a tiling of the plane,
using tiles from a tileset~$\mathcal{T}$ which is closed under
reflection, is said to be \emph{valid} if two adjacent tiles meet along an edge with the
same color and orientation and two tiles that are the reflection of
each other are never adjacent. In the following, such a tiling will be called of~\emph{Kari-Papasoglu type}.

If we consider only the horizontal
(vertical) symmetry we say that $\mathcal{T}$ is $h$-reflection-closed
(respectively, $v$-reflection-closed) tileset. Note that if
$\mathcal{T}$ is $h$-reflection-closed, then it is $ws$-deterministic
($es$-deterministic) if and only if it is $wn$-deterministic
(respectively, $en$-deterministic). Similarly, if $\mathcal{T}$ is
$v$-reflection-closed, then it is $ws$-deterministic
($wn$-deterministic) if and only if it is $es$-deterministic
(respectively, $en$-deterministic). Hence, if $\mathcal{T}$ is reflection-closed and $xy$-deterministic for some~$(x,y)\in\{(e,n),(e,s),(w,n),(w,s)\}$, then $\mathcal{T}$ is necessarily $4$-way deterministic. In~\cite{KaPa} the authors raised the problem of finding a $4$-way deterministic tileset which is valid, aperiodic and reflection-closed. Such a tileset would give an example of a CAT(0) complex whose fundamental group is not hyperbolic and does not contain a subgroup isomorphic to $\mathbb{Z}^{2}$, see~\cite{Gromov, KaPa}.\medskip
 
In this setting we can prove a statement analogous to Proposition \ref{prop: equivalence positive continuity}. We will prove that the search for aperiodic
$h$-reflection-closed tilesets that are $ws$- and $wn$-deterministic
is related to the search for 
automata in the class~$\SinkA$ whose set~$\krit$ of singular points is empty
and which are elementary-free, in the following sense. \\ 

\newcommand\hMaxSync{\small\textbf{{h-Maximal~non~Synchronizing }}}
\newcommand\KPMaxSync{\small\textbf{{4-way Maximal~non~Synchronizing }}}

First we need an analogous to Proposition~\ref{prop: no tiling characterization} for $h$-reflexion-closed tilings (resp. 4-way deterministic tilings). Let us define the \emph{prefix-reduced relation}: we say that $u \in \wt{\QQ}$ is smaller than $u'\in \wt{\QQ}$ for the prefix-reduced relation if $u \leq_p u'$ holds and both $u$ and~$u'$ are reduced.
 We define a property that will serve for the characterization of Kari-Papasoglu type tilings:
an automaton~$\auta=(\QQ,\XX,\xcdot,\xcirc) \in \SinkA $ satisfies property {\hMaxSync} (resp. {\KPMaxSync}) if there exists an integer $m$ such that, for all $v \in \XX^{m}$ (resp. reduced $v \in \wt{\XX^{m}}$), any  reduced $u \in \wt{\QQ^*}$  maximum in $\NSyn_{\wt{\QQ}}(v)$ for the prefix-reduced relation satisfies $uq \xcirc v\in \Syn(\auta)$ for all~$q\in \wt{\QQ}$.
\begin{proposition}\label{prop: no tiling characterization }
Let $\auta=(\QQ,\XX,\xcdot,\xcirc) \in \SinkA$. Then $\overline{\mathcal{T}}$ admits an $h$-reflexion-closed tiling (4-way deterministic tiling) if and only if property~{\hMaxSync} (resp.~{\KPMaxSync}) does not hold.
\end{proposition}
\begin{proof}
Similar to the proof of Proposition~\ref{prop: no tiling characterization}, but avoiding patterns~$xx^{-1}$ that are not allowed in the Kari-Papasoglu type tilings.
\end{proof}

A group generated by an automaton in the class~$\SinkA$ is said to be {\it elementary-free} if the only relations that it contains are words whose restrictions become eventually all trivial, \emph{i.e.} the set of its relations  may be described as the set of words $u\in \wt{(\QQ\smallsetminus\{e\})}^{*}$ with $\pi(u)=\id$, such that there exists an $n\ge 1$ for which $\overline{u\xcdot v}=e^{k_v}$ for every $v\in\XX^{\ge n}$ and some integer $k_v$ such that $|k_v|\leq |v|$. 
 
In this context, we say that a helix graph
$\wt{\mathcal{H}}_{k,n}$ of an automaton~$\auta\in\SinkA$
is \emph{strongly-singular} whenever any cycle~$(u,v)\mapright{}(u',v')\mapright{}\cdots\mapright{}(u,v)$
implies either~$u\in \wt{\{e\}}^{*}$ or~$\pi(u)=\id$.
We have the following proposition analogous to~Proposition~\ref{prop: equivalence positive continuity}.

\begin{proposition}\label{prop: equivalence continuity}
Let $\auta\in\SinkA$. The following are equivalent:
\begin{enumerate}[label=(\roman{enumi})]
\item \label{equ-cont1}$\gauta$ is elementary-free and the set~$\krit$ of singular points is empty;
\item \label{equ-cont2}$\gauta$ is elementary-free and for any $v\in \XX^{*}$, $g\in\St_{\gauta}(v^{\omega})$,
there is $n\ge 1$ with $g\xcdot v^{\omega}[:n]=\id$;
\item \label{equ-cont3}for any~$k,n\ge 1$, the helix graph~$\wt{\mathcal{H}}_{k,n}$ is strongly-singular;
\item \label{equ-cont4}there is no non-elementary pair $u\in \wt{(\QQ\smallsetminus\{e\})}^{*}, v\in \XX^{*}$ of commuting words.
\end{enumerate}
\end{proposition}

\begin{proof}
\ref{equ-cont1}Equivalence~$\Leftrightarrow$\ref{equ-cont2} follows from Lemma~\ref{lemmarestric}. Equivalence~\ref{equ-cont2}$\Leftrightarrow$\ref{equ-cont3} may be proven in an analogous way as in Proposition~\ref{prop: equivalence positive continuity}. Equivalence~\ref{equ-cont1}$\Leftrightarrow$\ref{equ-cont4} is a consequence of Proposition~\ref{prop: not continuous with inverses}.
\end{proof}

Note that $\gauta$ also acts naturally on $\wt{\XX}^*$, in what follows we consider this action.
Following \cite{DaRo15}, we say that a point $\xi \in \XX^{\omega}$ is \emph{essentially non-trivial} when $\lim|\overline{\xi[:n]}|\to_{n\to \infty} +\infty$. Moreover we say that a helix graph is essentially-singular helix whenever, if $(u,v)\mapright{}(u',v')\mapright{}\cdots\mapright{}(u,v)$ is a cycle, then $v^{\omega}$ is essentially non-trivial and either $u\in \wt{\{e\}}^{*}$, or $\pi(u)=\id$. We have the following proposition analogous to Proposition~\ref{prop: equivalence positive continuity}:

\begin{proposition}\label{prop: equivalence continuity 4way}
Let $\auta\in\SinkA$. The following are equivalent:
\begin{enumerate}[label=(\roman{enumi})]
\item \label{equ-cont1}$\gauta$ is elementary-free and the set of essentially non-trivial singular points is empty;
\item \label{equ-cont2}$\gauta$ is elementary-free and for any $v\in \wt{\XX^{*}}$ such that $v^{\omega}$ is essentially non-trivial, $g\in\St_{\gauta}(v^{\omega})$, there is $n\ge 1$ with $g\xcdot v^{\omega}[:n]=\id$;
\item \label{equ-cont3}for any~$k,n\ge 1$, the helix graph~$\wt{\mathcal{H}}_{k,n}$ is essentially singular;
\item \label{equ-cont4}there is no non-elementary pair $u\in \wt{(\QQ\smallsetminus\{e\})}^{*}, v\in \wt{\XX^{*}}$ of commuting words.
\end{enumerate}
\end{proposition}

\begin{proof}
Similar to the proof of Prop.~\ref{prop: equivalence continuity}.
\end{proof}

Putting together all the previous results we may characterize aperiodic $h$-reflection-closed tilesets (resp.  $h$-reflection-closed tilesets).  
\begin{theorem}\label{theo: ws-wn tilings}
With the above notation. The following are equivalent:
\begin{enumerate}[label=(\roman{enumi})]
\item \label{ws-wn1} there is a $ws$- and $wn$-deterministic tileset~$\mathcal{T}$ which is $h$-reflection-closed (resp. a 4-way deterministic which is  $h$- and $v$-reflection-closed)  that  tiles the discrete plane with aperiodic tilings of Kari-Papasoglu type;
\item \label{ws-wn2} there is an %invertible
automaton $\auta\in\SinkA$ such that $\gauta$ is elementary-free, the set of singular points (resp. essentially non-trivial singular points) is empty, and, in the automaton  $\mathrsfs{B}$  obtained from~$\auta\sqcup \auta^{-1}$ identifying the two sinks $e$ and $e^{-1}$, property~{\hMaxSync} (resp. {\KPMaxSync}) does not hold.
\end{enumerate}
\end{theorem}

\begin{proof}
If $\mathcal{T}$ is $h$-reflection-closed, we may fix a direction and divide the colors of the vertical edges of
$\mathcal{T}$ into two distinct and disjoint sets $\QQ\sqcup
\QQ^{-1}$, while we put for~$\XX$ the set of colors of the horizontal
edges of the tiles in $\mathcal{T}$. Note that, for each tile~$(q,a,p,b)$,
$q,p\in \wt{\QQ}$, $a,b\in \XX$, the corresponding horizontally
reflected tile is~$(q^{-1}, a, p^{-1}, b)$. The partition~$\QQ\sqcup
\QQ^{-1}$ induces a partition~$\mathcal{T}^{+}\sqcup\mathcal{T}^{-}$ on~$\mathcal{T}$
in the obvious way. Consider an associated %invertible
automaton $\auta\in\SinkA$ such that
$\overline{\mathcal{T}}(\auta)=\mathcal{T}^{+}$ (as in Section
\ref{sec: tiling comm}). Note that
$\overline{\mathcal{T}}(\auta\sqcup
\auta^{-1})=\mathcal{T}$.

Conversely, to any %invertible
automaton~$\auta\in\SinkA$, the tileset
$\overline{\mathcal{T}}(\auta\sqcup \auta^{-1})$ is
$ws$- and $wn$-deterministic, and it is $h$-reflection-closed. By an
argument very similar to the proof of Propositions~\ref{prop: periodic tiling} and~\ref{prop: reduced tileset}, it is not difficult to see that,
in the previous correspondence, $\auta$ has a non-elementary pair~$u\in \wt{(\QQ\smallsetminus\{e\})}^{*}, v\in \XX^{*}$ of commuting words where $u$ is non-trivial and reduced if and only if the corresponding tileset~$\overline{\mathcal{T}}(\auta\sqcup \auta^{-1})$ admits a periodic tiling (in the sense of Kari-Papasoglu). Hence, the equivalence in the statement follows from Proposition~\ref{prop: no tiling characterization } for the existence of a tiling of Kari-Papasoglu type and Proposition~\ref{prop: equivalence continuity} for its aperiodicity. The proof for the 4-way case is similar and uses Proposition~\ref{prop: equivalence continuity 4way}.
\end{proof}

The last theorem gives a characterization of specific Wang tilings in the language of Mealy automata. This is another motivation to further explore this connection.

%-------------------------------------------------------------------------------------------------------------------------------------------------
%-------------------------------------------------------------------------------------------------------------------------------------------------
%-------------------------------------------------------------------------------------------------------------------------------------------------
\section{Some open problems}
\begin{prob}
Let $\auta\in\SinkA$ be an automaton generating a free group. Is it always the case that there is a point in the boundary whose Schreier graph is infinite?
\end{prob}

\begin{prob}\label{prob:isolated}
Given that $\auta$ is minimized, can two singular points have isomorphic infinite Schreier Graphs?
\end{prob}

\begin{prob}
Given a Mealy automaton~$\auta$, is it decidable whether $\mathcal{P}(\auta)=\varnothing$?
\end{prob}

\begin{prob}
Are there interesting classes of automata where the non-elementary commuting pair is decidable?
\end{prob}

\bibliographystyle{plain}
%\bibliography{schreier}
%\end{document}

\def\cprime{$'$}

\end{document}